\RequirePackage{fix-cm}
\documentclass{svjour3}   
\smartqed

\def\Frac#1#2{\frac{\displaystyle{#1}}{\displaystyle{#2}}}
\usepackage{graphicx}
\usepackage{amssymb}
\usepackage{amsfonts}

\begin{document}

\title{Fast and reliable high accuracy computation of Gauss--Jacobi quadrature\thanks{
This work was supported by Ministerio de Ciencia, Innovaci\'on y Universidades, Spain, 
projects MTM2015-67142-P (MINECO/FEDER, UE) and PGC2018-098279-B-I00 (MCIU/AEI/FEDER, UE).
}}

\titlerunning{Iterative computation of Gauss--Jacobi quadrature} 

\author{Amparo Gil \and Javier Segura \and Nico~M.~Temme}
\authorrunning{A. Gil, J. Segura, N. M.~Temme} 

\institute{Amparo Gil \at
              Departamento de Matem\'atica Aplicada y CC. de la Computaci\'on.
ETSI Caminos. Universidad de Cantabria. 39005-Santander, Spain \\
              \email{gila@unican.es}         
           \and
           Javier Segura \at
           Departamento de Matem\'aticas, Estad\'{\i}stica y Computaci\'on.
Facultad de Ciencias. Universidad de Cantabria. 39005-Santander, Spain.\\
\email{segurajj@unican.es} 
\and
Nico M.~Temme \at
IAA, 1825 BD 25, Alkmaar, The Netherlands. Former address: Centrum Wiskunde \& Informatica (CWI), 
        Science Park 123, 1098 XG Amsterdam,  The Netherlands.\\
\email{nicot@cwi.nl} 
}

\date{Received: date / Accepted: date}

\maketitle

\begin{abstract}
Iterative methods with certified convergence for the computation of Gauss--Jacobi quadratures 
are described. The methods do not require a priori estimations of the nodes to guarantee its fourth-order
convergence. They are shown to be generally faster than previous methods and without practical restrictions on
the range of the parameters. The evaluation of the nodes and weights of the quadrature is
exclusively based on convergent processes which, together with the fourth order convergence of the fixed
point method for computing the nodes, makes this an ideal approach for high accuracy computations, so much
so that computations of quadrature rules with even millions of nodes and thousands of digits are possible in 
a typical laptop. 
\end{abstract}


\section{Introduction}

Given an integral $I(f)=\int_{a}^b f(x) w(x) dx$, with $w(x)$ a weight function in the interval $[a,b]$, it is
said that the $n$-point quadrature rule $Q(f)=\sum_{i=1}^n w_i f(x_i)$ is a Gaussian quadrature if
it has the highest possible degree of exactness, that is, if $I(f)=Q(f)$ for all polynomials of degree smaller
than $2n$.

Gauss--Jacobi quadrature is, together with Gauss--Hermite and Gauss--Laguerre quadratures, 
one of the three classical Gauss
quadrature rules and it is, without any doubt, the most widely used of them. This rule corresponds
to the weight function $w(x)=(1-x)^{\alpha} (1+x)^{\beta}$, $\alpha,\beta>-1$, in the interval $[-1,1]$ and
it has as particular cases Gauss-Chebyshev quadratures ($|\alpha|=|\beta|=1/2$) and Gauss--Legendre quadrature
$\alpha=\beta=0$.

Because of the optimal degree of exactness of the Gauss rules, they have fast convergence as the degree 
increases (specially for analytic functions), and they are one of the most popular methods of numerical integration,
appearing in countless applications. However, Gauss rules are usually seen as hard to compute and for this
reason alternative simpler rules as the Clenshaw-Curtis rules may be preferred 
\cite{Tre:2008:IGQ}. Nevertheless, Gauss rules are in any case
optimal in terms of degree of exactness, and for integrals where the explicit weight of the quadrature appears 
(like $w(x)=(1-x)^{\alpha}(1+x)^{\beta}$
for Gauss--Jacobi) they are difficult to beat. In addition, the efficiency of the computation of Gauss rules has
 dramatically improved in recent years as shown for instance in  
\cite{Glaser:2007:AFA,Hale:2013:FAA,Bogaert:2014:IFC,Gil:2018:GHL,Gil:2019:NIC,Gil:2019:FRA,Bremer:2019:FAJ}.
Let us first briefly describe these methods, in particular
for Gauss--Jacobi quadrature: 

1. The Golub--Welsch algorithm \cite{Golub:1969:COG} is a simple approach based on the diagonalization 
of the Jacobi matrix associated to the three-term recurrence relation. For small degree $n$ 
it is viable method, but its complexity scales as ${\cal O}(n^2)$ and it becomes slow as the degree increases.

2. Iterative methods: they 
use the fact that the nodes of Gaussian quadrature are the roots of the 
orthogonal polynomials associated to the quadrature, while the weights are related to the derivative of the 
polynomial at the nodes. The complexity increases linearly. 
This is the approach considered in \cite{Hale:2013:FAA}, where the associated orthogonal polynomials 
(Jacobi polynomials) are computed by means of asymptotic formulas for large $n$, as well as the initial values for the nodes for starting the Newton iteration. The main limitation of this approach seems to be that the initial values 
for starting the Newton
iterations only guarantee convergence for small $|\alpha|$ and $|\beta|$. 
For the particular case of Gauss--Legendre quadrature ($\alpha=\beta=0$) there is a larger number
 of iterative methods available, see 
\cite{Yakimiw:1996:ACW,Swarztrauber:2002:OCT,Glaser:2007:AFA,Bogaert:2012:COL,JOH:2018:FAR}.  
A recent alternative iterative method for Gauss--Jacobi quadrature is that of \cite{Bremer:2019:FAJ}, 
which also appears to be faster than the approach of \cite{Hale:2013:FAA}, but it is more limited than \cite{Hale:2013:FAA} with respect to the parameters 
($|\alpha|,|\beta|<1/2$).

3. Asymptotic methods: explicit approximations for the Gauss--Jacobi nodes and weights which do not require 
iterative refinement are given in \cite{Gil:2019:NIC}. These methods are faster than iterative methods
precisely because no iterations are needed and explicit formulas are used instead. 
As for the case of iterative methods, these approximations have limited
validity and other type of asymptotic formulas should be considered for large $\alpha$ and/or $\beta$; 
a first
step in this direction is given in \cite{Gil:2019:AEO}. A previous asymptotic method for the particular case of
Gauss--Legendre quadrature is given in \cite{Bogaert:2014:IFC}. For an analysis of asymptotic methods for generalized 
Gauss quadratures (Gauss-Jacobi in particular) based on Riemann-Hilbert analysis see \cite{Opsomer:2018:AFO}.

4. Global high-order iterative methods with certified convergence: 
for Gauss--Hermite and Gauss--Laguerre quadratures, a new approach was considered in \cite{Gil:2019:FRA} 
which combines the use 
of the fourth-order globally convergent fixed point method of \cite{Segura:2010:RCO} 
with the use of local Taylor series for computing the weights 
(Taylor series are also considered in \cite{Glaser:2007:AFA}). This approach 
produced a fast, reliable and unrestricted algorithm which outperformed previous methods 
in terms of speed (though the asymptotic methods in \cite{Gil:2019:FRA} may be faster for large degrees), 
accuracy and available range of computation.

In this paper, we complete the construction of fast methods for classical Gauss quadratures with the description
of high-order iterative methods for Gauss--Jacobi quadrature. We therefore close the analysis of 
classical quadratures, adding to the asymptotic methods \cite{Gil:2018:GHL} (Gauss--Hermite and Gauss--Laguerre) 
and \cite{Gil:2019:NIC} (Gauss--Jacobi), 
and to the high-order iterative method of \cite{Gil:2019:FRA} (Gauss--Hermite and Gauss--Laguerre), 
the corresponding high-order iterative method for Gauss--Jacobi. 
We believe that an optimal algorithm for the computation of Gauss quadratures
in fixed precision will involve both the asymptotics-free iterative methods and the
iteration-free asymptotic methods which are completed in the present paper.

As advanced in \cite{Gil:2019:FRA}, the implementation
of the global iterative methods for Gauss--Jacobi 
is not so straightforward as for the Hermite and Laguerre cases, not only because there are more 
parameters involved, but also because the possible changes of variable for the 
 Liouville transformations of the ODE needed in the method 
are not amenable to the use of Taylor series. In 
practical terms this, as we will see, means that we will have to combine different fixed point methods associated
with different Liouville transformations and an independent application of Taylor series.
As we will see, our methods are generally faster than \cite{Hale:2013:FAA} and 
with a much larger range of validity than 
\cite{Hale:2013:FAA,Bremer:2019:FAJ}, and they have no rival for high accuracy computations due to his
high order of convergence. It is
also a much simpler method than previous methods, particularly for the symmetric case 
(that is, for Gauss--Gegenbauer quadrature).

The structure of the paper is as follows. Firstly, we describe the main ingredients of the method and summarize 
the relations satisfied by Jacobi polynomials that will be used in the paper. In the second place, we describe
a basic algorithm for the symmetric case $\alpha=\beta\ge 0$ (Gauss--Gegenbauer, including Gauss--Legendre). 
For this simple symmetric case some accuracy problems, however, appear when the parameters $\alpha$ and/or $\beta$ are close
to $-1$, which require further attention. Next, we describe the more general algorithm for 
Gauss--Jacobi quadrature for $\alpha,\beta >-1$, which adds two features with respect to the symmetric case: a starting 
procedure based on the three-term recurrence relation and an alternative Liouville transformation for the 
extreme nodes for negative parameters; this modification solves the numerical accuracy problems for parameters $\alpha$ and
$\beta$ approaching $-1$.
Finally, we provide numerical evidence of the speed and accuracy of the method, including very high precision 
computations for the symmetric case (even with more than 1000 digits). 
We compare our method
against the chebfun \cite{Dris:2014:CG} 
implementation of the methods in \cite{Hale:2013:FAA,Bogaert:2014:IFC} in the regions of parameters where those 
are valid. As we will discuss, our method is competitive in speed and accuracy with previous methods (and notably faster 
for the cases $\alpha\neq \beta$) and with the advantage that it works without practical restrictions on the
parameters. It has the additional benefit that computations  with very high accuracy are possible and they can be efficiently
 performed thanks to the high--order convergence of the method and the fact that is is based on convergent processes, 
different to the asymptotic approaches of \cite{Bogaert:2014:IFC,Hale:2013:FAA,Gil:2019:NIC}.

\section{Basic ideas and main formulas}

Our algorithm, as all the other iterative methods, is based on two well-known facts. The first one is that
the nodes of the Gauss--Jacobi quadrature of degree $n$ are the roots $x_i$, $i=1,\ldots n$, 
of the Jacobi polynomial of degree $n$,
$P_n^{(\alpha,\beta)}(x)$, and the second is that the weights can be computed in terms of the derivative 
at the nodes as
\begin{equation}
\label{peso}
w_i=\Frac{M_{n,\alpha,\beta}}{(1-x_i^2)(P_n^{(\alpha,\beta) '}(x_i))^2},
\end{equation}
with
\begin{equation}
\label{constante}
M_{n,\alpha,\beta}=2^{\alpha+\beta+1}\Frac{\Gamma (n+\alpha+1)\Gamma (n+\beta+1)}{n!\Gamma (n+\alpha+\beta+1)}.
\end{equation}
The Jacobi polynomials can be written in terms of Gauss hypergeometric functions as
\begin{equation}
\label{defh}
P_n^{(\alpha,\beta)}(x)=\Frac{(\alpha+1)_n}{n!}\,_{2}{\rm F}_{1}
\left(-n,n+\alpha+\beta+1;\alpha+1;\Frac{1-x}{2}\right),
\end{equation}
and they satisfy the symmetry relation
\begin{equation}
\label{sym}
P_n^{(\alpha,\beta)}(x)=P_n^{(\beta,\alpha)}(-x).
\end{equation}

For computing the nodes, we use the global fixed point method of \cite{Segura:2010:RCO} which applies to second
 order homogeneous linear ODEs. As it is well known, we have 
\begin{equation}
\label{odeJ}
(1-x^2)y''(x)+ \left[(\beta-\alpha)-(2+\alpha+\beta)x\right]y'(x)+n(n+\alpha+\beta +1)y(x)=0
\end{equation}
for $y=P_n^{(\alpha,\beta)}(x)$. Starting from this equation, we will consider several Liouville transformations 
which lead to equations in normal form, suitable for applying the fixed point method \cite{Segura:2010:RCO}.

The main results of \cite{Segura:2010:RCO} that we will use in our methods can be condensed
in the following theorem (where dots mean
derivative with respect to $z$):
\begin{theorem}
\label{fixed}
Let $Y(z)$ be a solution of $\ddot{Y}(z)+\Omega (z)Y(z)=0$ and let ${\protect\rm a}$ be 
such that $Y({\protect\rm a} )=0$. 
Let ${\protect\rm b}\neq {\protect\rm a}$ such that $\Omega ({\protect\rm b})>\Omega ({\protect\rm a})$ and 
$Y(z)\neq 0$ in the open interval $I$ between ${\protect\rm a}$ and ${\protect\rm b}$. 
Assume that $\Omega(z)$ is differentiable and monotonic in the closure of $I$. Let 
$j=\mbox{sign}({\protect\rm b}-{\protect\rm a})$, then for any $z^{(0)}\in I\bigcup \{{\protect\rm b}\}$, 
the sequence $z^{(i+1)}=T_j(z^{(i)})$, $i=0,1,\ldots$, with
\begin{equation}
\label{fixedpoint}
T_{j}(z)=z-\Frac{1}{\sqrt{\Omega (z)}}\arctan_{j}\left(\sqrt{\Omega (z)}\Frac{Y(z)}{\dot{Y}(z)}\right)
\end{equation}
and
\begin{equation}
\label{arctan}
\arctan_{j}(\zeta)=\left\{
\begin{array}{l}
\arctan(\zeta)
\mbox{ if }j\zeta > 0,\\ 
\arctan(\zeta)+j\pi \mbox{ if }j\zeta\le 0,\\ 
j\pi/2 \mbox{ if } \zeta=\pm \infty,
\end{array}
\right.
\end{equation}
is such that $\{z^{(i)}\}_{i=1}^{\infty}\subset I$ and it converges monotonically to the root 
${\protect\rm a}$ with order
of convergence $4$ and asymptotic error constant $\dot{\Omega}({\protect\rm a})/12$, that is:
$$
\displaystyle\lim_{i\rightarrow \infty}\Frac{z^{(i+1)}-{\protect\rm a}}{(z^{(i)}-{\protect\rm a})^4}
=\Frac{\dot{\Omega}({\protect{\rm a}})}{12}.
$$

\end{theorem}

\begin{remark}{\it
Observe that in the previous theorem ${\protect\rm b}$ could be such that $Y({\protect\rm b})=0$. 
Therefore the theorem gives a procedure
to compute zeros in succession in the direction of decreasing values of $\Omega (z)$. If $Y(z^{(0)})=0$ then 
the first iteration is $z^{(1)}=z^{(0)}-j\pi/\sqrt{\Omega (z^{(0)})}$.}
\end{remark}

\begin{remark}{\it
Because the fixed point method generates monotonic sequences when $\Omega (z)$ is monotonic, it does not 
show local convergence around each zero, but only lateral convergence (which is why the previous remark is
true). This means that if, for instance, $\Omega (z)$ is decreasing in an interval and 
$a_0$ and $a_1$ are two consecutive roots in the interval, 
$a_0<a_1$, then the iteration converges to $a_0$ for values of $z$ close enough 
to $a_0$ and such that $z<a_0$, but it will converge to $a_1$ for
starting values in $[a_0,a_1)$. 
If instead of this, a method with bilateral local convergence
 is needed, we only need to replace the definition of (\ref{arctan}) by the usual definition of the arctangent;
that is, one can consider the fixed point method
\begin{equation}
\label{fixedpoint2}
g(z)=z-\Frac{1}{\sqrt{\Omega (z)}}\arctan \left(\sqrt{\Omega (z)}\Frac{Y(z)}{\dot{Y}(z)}\right).
\end{equation}
This redefinition of the 
fixed point method converges in wide intervals around each root under mild assumptions (see Theorem 
3.2 of  \cite{Segura:2010:RCO}).
}
\end{remark}

Theorem \ref{fixed} is the main tool for computing the nodes of Gauss--Jacobi quadrature. For applying this
result we first need to transform our equation to normal form, suppressing the first derivative term by 
means of a Liouville transformation; in addition, we will need a method to compute $y(z)/\dot{y}(z)$. 
We summarize next the Liouville transformations used in our algorithms, and later we discuss the methods
of computation.

\subsection{Three Liouville transformations of the Jacobi equation}

Let
$$
\label{ode}
y''(x)+B(x)y'(x)+A(x)y(x)=0.
$$
We consider a change of variables followed by a transformation to normal form so that the 
transformed equation reads
\begin{equation}
\label{norm}
\ddot{Y}(z)+\Omega(z)Y(z)=0,
\end{equation}
where dots represent the derivative with respect to $z$. In terms of the original variable $x$ 
(see for instance \cite{Dea:2004:NIF}) we have
$$
Y(z(x))=\sqrt{z'(x)}\exp\left(\frac12\int B(x)dx\right)y(x)
$$
and
$$
\Omega (z(x))=\Frac{1}{z'(x)^2}\left(A(x)-\Frac{B'(x)}{2}-\Frac{B(x)^2}{4}+\Frac{3z''(x)^2}{4z'(x)^2}
-\Frac{z'''(x)}{2 z'(x)}\right).
$$

For the Jacobi equation (\ref{odeJ}) we have $B(x)=\Frac{\beta+1}{1+x}-\Frac{\alpha+1}{1-x}$ and then
$$
Y(z(x))=\sqrt{z'(x)}(1-x)^{(1+\alpha) /2} 
(1+x)^{(1+\beta)/2} y(x).
$$

Because, according to Theorem \ref{fixed}, the monotonicity properties of $\Omega (z)$ are needed 
in order to apply the fixed point method, the changes of variable to be considered should allow a
simple determination of these properties. In \cite{Dea:2004:NIF,Dea:2007:GSI} the changes of variable for which
the determination of these properties reduces to the solution of a second order algebraic equation 
are analyzed systematically. Of these, we will use the three symmetric changes of variable in terms of 
elementary functions described in \cite{Dea:2004:NIF}, which are those with $z'(x)=(1-x^2)^{p}$ with
$p=0,-1/2,-1$.

\subsubsection{Trivial transformation}
\label{trivi}

For $p=0$ we have the trivial change $z(x)=x$. For later convenience we denote the transformed function
as $\tilde{Y}$ instead of $Y$. The transformed function in this case is 
\begin{equation}
\label{Ytil}
\tilde{Y}(x)=(1-x)^{(\alpha+1)/2}(1+x)^{(\beta+1)/2} P_n^{(\alpha,\beta)}(x)
\end{equation}
and satisfies
\begin{equation}
\label{odetr}
\begin{array}{c}
\tilde{Y}''(x)+\Omega(x)\tilde{Y}(x)=0,\\
\\
\Omega(x)
%= \Frac{1}{4(1-x^2)^2}\left[(1-L^2)x^2+2(\beta^2-\alpha^2)x+L^2-2\alpha^2-2\beta^2+3\right]
%\\
=\Frac{(L^2-1)(1-x^2)-2(\alpha^2 -1)(1+x)-2(\beta^2 -1)(1-x)}{4(1-x^2)^2},
\end{array}
\end{equation}
where $L=2n+\alpha+\beta+1$.

As described in \cite{Dea:2007:GSI}, the monotonicity properties of $\Omega (x)$ are not simple, and therefore
this transformation is of no use for the fixed point method. However, as we will see, it will be useful for 
applying Taylor series in order to compute the function $Y(z)/\dot{Y}(z)$ appearing in the fixed point method
 (\ref{fixed}).

Notice that in terms of the derivative of $\tilde{Y}$ at the nodes, the weights can be written as
\begin{equation}
\label{pesoz}
w_i=\Frac{M_{n,\alpha,\beta}}{\tilde{Y} '(x_i)^2} (1-x_i)^{\alpha}(1+x_i)^{\beta},
\end{equation}
and that, given that $\tilde{Y}$ satisfies an equation in normal form, the
quantity $\omega_i =1/\tilde{Y} '(x_i)^2$ is well conditioned as a function of the node 
$x_i$ because the function $\omega(x)=1/\tilde{Y}'(x)^2$ is such that $\omega'(x_i)=0$.
The main
source of error for the weights will be in the factor $(1-x_i)^{\alpha}(1+x_i)^{\beta}$,
particularly for the nodes close to $\pm 1$.

As we will discuss later, these transformation will be used for computing most of the nodes and weights, 
 and we will use Taylor series based on (\ref{odetr}). Consequently, the weights will be computed through 
(\ref{pesoz}), which uses the scaled weight. The method thus naturally computes the scaled weights and
from them the unscaled weights. The scaled weights are not only better conditioned that the unscaled weights,
 but they are also less prone to underflow for large values of $\alpha$ and $\beta$, which allows for a
 better control of these type of problems.

\subsubsection{Angular transformation}
\label{angular}

For $p=-1/2$ we denote the new variable by $\theta$ instead of $z$. We have $x=\cos \theta$, 
$\theta\in [0,\pi]$ and the transformed function 
\begin{equation}
\label{yxa}
Y(\theta(x))=(1-x)^{(\alpha+1/2)/2}(1+x)^{(\beta+1/2)/2} P_n^{(\alpha,\beta)}(x)
\end{equation}
satisfies equation (\ref{norm}), with $\theta\equiv z$, and
\begin{equation}
\label{eyx}
\begin{array}{ll}
\Omega(\theta(x))=\Frac{1}{4}L^2-\Frac{\alpha^2-1/4}{2(1-x)}-\Frac{\beta^2-1/4}{2(1+x)}.
\end{array}
\end{equation}

In terms of the transformed function $Y(\theta)$ the weights can be written as
$$
w_i =\Frac{M_{n,\alpha,\beta}}{\dot{Y}(\theta_i)^2}(1-x_i)^{\alpha+1/2}(1+x_i)^{\beta+1/2},
$$
where the dot means derivative with respect to $\theta$. In this expression 
$\omega_i=1/{\dot{Y}(\theta_i)^2}$ is well conditioned as a function of $\theta_i =\arccos(x_i)$ and has 
a slow variation as a function of the nodes $\theta_i$ as 
$n$ becomes large. Indeed, according to the circle theorem \cite{Davis:1961:SGT}
$$
w_i\sim \Frac{\pi}{n}w(x_i)\sqrt{1-x_i^2},\quad n\rightarrow \infty,
$$ 
with $w(x)$ the weight function, and for Gauss--Jacobi quadrature this gives
\begin{equation}
\label{circle}
w_i\sim \Frac{\pi}{n} (1-x_i)^{\alpha+1/2} (1+x_i)^{\beta+1/2}.
\end{equation}

The monotonicity properties of $\Omega (x)$ are simple to analyze: $\Omega (x)$ has one minimum in $(-1,1)$ when
$|\alpha|>1/2$ and $|\beta| >1/2$, one maximum when $|\alpha|<1/2$ and $|\beta|<1/2$ and it is monotonic
in the rest of cases. It is possible to construct methods for computing the Gauss--Jacobi quadratures by
using this transformation, however we will prefer the next transformation ($p=-1$) because the monotonicity
properties are even simpler. In some cases we will use this angular transformation for computing
the extreme nodes.

\subsubsection{Transformation to ${\mathbb R}$}
\label{erre}

With $p=-1$ we have the change $x=\tanh z$, $z\in {\mathbb R}$, and the transformed function 
\begin{equation}
\label{yxz}
Y(z(x))=(1-x)^{\alpha/2}(1+x)^{\beta/2} P_n^{(\alpha,\beta)}(x)
\end{equation}
satisfies equation (\ref{norm}) with
\begin{equation}
\label{Ozx}
\begin{array}{ll}
\Omega(z(x))
& =\Frac{1}{4}\left[(L^2-1)(1-x^2) -2 \alpha^2 (1+x)-2 \beta^2 (1-x))\right].
\end{array}
\end{equation}

In terms of the derivative with respect to $z$ at the nodes the weights can be written
$$
w_i =\Frac{M_{n,\alpha,\beta}}{\dot{Y}(z_i)^2}(1-x_i)^{\alpha+1}(1+x_i)^{\beta+1},
$$
where $z_i=\tanh^{-1} (x_i)$.

The coefficient $\Omega (z(x))$ 
has a maximum at $x_e=(\beta^2-\alpha^2)/(L^2-1)$, for any values of $\alpha$ and $\beta$. 
Because of these simple monotonicity properties, we will use this transformation for 
our method. According to Theorem 
\ref{fixed}, the fixed point method has to be applied in the direction of decreasing $\Omega (z)$.
Then the method can proceed starting
at $z_e= z(x_e)$, with a forward sweep for $z>z_e$ and a backward sweep for $z<z_e$. This is similar 
to the procedure for Gauss--Hermite quadrature, in particular for the case $\alpha=\beta$, when the
method starts at $x=0$ and the problem is symmetric.

\subsection{Methods of computation}

As basic method of 
computation our algorithms will use local Taylor series, however, alternative 
methods (recurrences, continued fraction) are employed
for the non-symmetrical case $\alpha \neq \beta$, 
and also for the extreme zeros for negative parameters. We start summarizing some 
information on the recurrences and later we describe the use of Taylor series.

\subsubsection{Recurrence relations and continued fractions}

Using (\ref{odeJ}) and the differentiation formula
$$
\frac{d}{dx}P_n^{(\alpha,\beta)}(x)=\Frac{n+\alpha+\beta+1}{2} P_{n-1}^{(\alpha+1,\beta+1)}(x),
$$
we obtain the recurrence relation
$$
\begin{array}{ll}
(n+\alpha+\beta+1)(1-x^2)P_{n-1}^{(\alpha+1,\beta+1)}(x)&+2[\beta(1-x)-\alpha(1+x)]
P_{n}^{(\alpha,\beta)}(x)\\
&+4(n+1)P_{n+1}^{(\alpha-1,\beta-1)}(x)=0,
\end{array}
$$
which can be used to compute $P_{n}^{(\alpha,\beta)}(x)$ starting from
$$
P_0^{(\alpha+n,\beta+n)}(x)=1,\quad P_1^{(\alpha+n-1,\beta+n-1)}(x)=\frac12\left(\alpha-\beta+(\alpha+\beta+2n)x\right),
$$
and is an alternative to the more popular three-term recurrence relation
\begin{equation}
\label{TTRR}
\begin{array}{l}
2(n+1)(n+\alpha+\beta+1)(2n+\alpha+\beta)P_{n+1}^{(\alpha,\beta)}(x)\\
-(2n+\alpha+\beta+1)\left[(2n+\alpha+\beta)(2n+\alpha+\beta+2)x+\alpha^2-\beta^2\right]P_n^{(\alpha,\beta)}(x)\\
+2(n+\alpha)(n+\beta)(2n+\alpha+\beta+2)P_{n-1}^{(\alpha,\beta)}(x)=0.
\end{array}
\end{equation}

Other relations can be found which are also useful for computing polynomial ratios, although they do not lead to finite exact recurrence methods for the
polynomials.  One of them is the recurrence relation
\begin{equation}
\label{reca}
\begin{array}{l}
(n+\alpha+\beta+1)(1-x) P_n^{(\alpha+1,\beta)}(x)\\-\left[(2n+\alpha+\beta+1)(1-x)+2\alpha\right]
P_n^{(\alpha,\beta)}(x)+2(\alpha+n) P_n^{(\alpha -1,\beta)}(x)=0,
\end{array}
\end{equation}
which leads to a continued fraction that will be useful in our algorithms.
This recurrence corresponds to the case of the $(0\,+\,+)$ for hypergeometric functions 
(using (\ref{defh})), and from the analysis of this recurrence (see \cite{Gil:2007:NSS})
we deduce that $P_n^{(\alpha,\beta)}(x)$ is minimal as $\alpha\rightarrow \infty$ in the disc in the
complex plane $|x-1|<2$, and in particular for $x\in (-1,1)$. This, on account of Pincherle's theorem 
\cite[Thm. 4.7]{Gil:2007:NMF},
means that the ratio $P_n^{(\alpha,\beta)}(x)/P_n^{(\alpha-1,\beta)}(x)$ can be computed via a 
continued fraction. By re-writing the previous recurrence as
$$
H_{\alpha}=\Frac{a_{\alpha}}{b_{\alpha}+H_{\alpha+1}},
$$
where $H_{\alpha}=P_n^{(\alpha,\beta)}(x)/P_n^{(\alpha-1,\beta)}(x)$ and
\begin{equation}
a_{\alpha}=-\Frac{2(\alpha+n)}{(n+\alpha+\beta+1)(1-x)},\quad b_{\alpha}=-1-\Frac{n(1-x)+
2\alpha}{(n+\alpha+\beta+1)(1-x)}.
\end{equation}
Iterating with have
\begin{equation}
\label{cfP}
\Frac{P_n^{(\alpha,\beta)}(x)}{P_n^{(\alpha-1,\beta)}(x)}=H_{\alpha}=\Frac{a_{\alpha}}{b_{\alpha}+}
\Frac{a_{\alpha+1}}{b_{\alpha+1}+}\cdots,
\end{equation}
which converges in $|x-1|<2$, and with faster convergence as we are closer to $x=1$. For this reason, it
will be an interesting method when computing the extreme zeros (notice that, because of (\ref{sym}), this can
also be used for nodes close to $x=-1$).

In connection with this recurrence relation over $\alpha$, we have the following relation for the derivative,
which will be used later
\begin{equation}
\label{dera}
(1-x^2)\frac{d}{dx}P_n^{(\alpha,\beta)}(x)=(n(1-x)+2\alpha)P_n^{(\alpha,\beta)}(x)
-2 (\alpha +n) P_n^{(\alpha -1,\beta)}(x).
\end{equation}

\subsubsection{Local Taylor series}
\label{taygen}

We start from the transformed ODE (\ref{odetr}), which we write as

$$
Q(x)\tilde{Y}''(x)+R(x)\tilde{Y}(x)=0,
$$
where, as before, $\tilde{Y}(x)=(1-x)^{(\alpha +1)/2} (1+x)^{(\beta +1)/2} P_n^{(\alpha,\beta)}(x)$,
and $Q(x)$ and $R(x)$ are the polynomials
\begin{equation}
\begin{array}{l}
Q(x)=4(1-x^2)^2,\\
\\
R(x)=(L^2-1)(1-x^2)-2(\alpha^2 -1)(1+x)-2(\beta^2 -1)(1-x).
\end{array}
\end{equation}

Given the initial values $Y(x)$ and $Y'(x)$, we can use Taylor series to compute the function and the derivative
at a different point $x+h$. 
The Taylor series centered at $x$ are
$$
\tilde{Y}(x+h)=\displaystyle\sum_{j=0}^{\infty}\Frac{u_j(x)}{j!}h^j,\,
\tilde{Y}'(x+h)=\displaystyle\sum_{j=0}^{\infty}\Frac{u_{j+1}(x)}{j!}h^j,
$$
provided $x+h$ is inside the interval of convergence around $x$. Of course, the series will be truncated
to a finite number of terms $N$. For computing these series, we need to evaluate the successive derivatives of 
$Y(x)$ starting from $Y(x)$ and $Y'(x)$.
For this purpose, we can differentiate the ODE. 

Differentiating $m$ times we have
$$
\displaystyle\sum_{k=0}^{4}
\left(
\begin{array}{c}
m\\
k
\end{array}
\right)
Q^{(k)} \tilde{Y}^{(2+m-k)}
+
\displaystyle\sum_{k=0}^{2}
\left(
\begin{array}{c}
m\\
k
\end{array}
\right)
R^{(k)} \tilde{Y}^{(m-k)}=0.
$$

This gives, denoting $u_k =\tilde{Y}^{(k)}$

\begin{equation}
\label{der1}
\begin{array}{l}
Q u_{j+2}+j Q' u_{j+1}+\left(\Frac{j(j-1)}{2}Q''+R\right)u_j\\
+ j\left(
\Frac{(j-1)(j-2)}{6}Q'''+R'\right)u_{j-1}\\
+ \Frac{j (j-1)}{2}\left(
\Frac{(j-2)(j-3)}{12}Q^{(4)}+R''\right)u_{j-2}=0.
\end{array}
\end{equation}

Instead of using the recurrence (\ref{der1}), in our algorithms we prefer to compute 
the quantities $a_j=u_j/j!$, which are less prone to overflow. The 
truncated Taylor series are then
\begin{equation}
\label{truncs}
\tilde{Y}(x+h)\simeq\displaystyle\sum_{j=0}^{N} a_j(x) h^j,\quad 
\tilde{Y}'(x+h)\simeq\displaystyle\sum_{j=0}^{N}(j+1) a_{j+1}(x) h^j,
\end{equation}
and the coefficients $a_j$ satisfy:
\begin{equation}
\label{recun}
\begin{array}{l}
(j+2)(j+1)Q a_{j+2}+(j+1) j Q' a_{j+1}+\left(\Frac{j(j-1)}{2}Q''+R\right)a_j\\
+ \left(
\Frac{(j-1)(j-2)}{6}
Q'''+R'\right)a_{j-1}\\
+ \Frac{1}{2}\left(
\Frac{(j-2)(j-3)}{12}
Q^{(4)}+R''\right)a_{j-2}=0,\,j=0,1,\ldots
\end{array}
\end{equation}
with $a_{-2}=a_{-1}=0$, $a_0=\tilde{Y}(x)$, $a_2=\tilde{Y}'(x)$. It is important to take into account that
the use of recurrence relations may be extremely unstable when the conditioning is not appropriate. It may
happen that there is a solution of the same recurrence $\{b_j\}$ such that the comparison with the wanted solution 
$\{a_j\}$ gives $\limsup_{n\rightarrow \infty}\sqrt[n]{|a_n/b_n|}<1$, in which case $\{b_n\}$ would dominate the forward application
of the recurrence and it may ruin the numerical computation, particularly if many terms of the series
are needed. However, it is  possible to prove that this can not happen in
most occasions (see Appendix A). 
This, although does not prove stability, at least disproves catastrophic exponential degradation of accuracy. 
Furthermore, the steps in the Taylor series will not be large and the number of terms required is not too large.
Numerical experiments indeed prove that the computation is stable.

\section{Gauss--Gegenbauer quadrature}
\label{GGQ}

We take now $\alpha=\beta=\lambda$ and consider the Liouville transformation 
with $x=\tanh z$. Then $Y(z)=\cosh (z)^{-\lambda}
P_n^{(\lambda,\lambda)}(\mbox{tanh}(z))$ satisfies 
\begin{equation}
\label{Yz}
\ddot{Y}(z)+\Omega (z)Y(z)=0,\,
\Omega(z)=\Frac{1}{4}\left[\Frac{L^2-1}{\mbox{cosh}^2 (z)}-4\lambda^2\right],
\end{equation}
and we can compute the roots of $Y(z)$, similarly as we did for Gauss--Hermite 
starting from $z=0$ and evaluating the positive roots in increasing order. By symmetry, the negative
roots are the same as the positive but with opposite sign.

The fixed point method (Theorem \ref{fixed}) is
\begin{equation}
\label{fixedz}
T(z)=z-\Frac{1}{\sqrt{\Omega (z)}}\arctan_{-1}\left(\sqrt{\Omega (z)} Y(z)/\dot{Y}(z)\right),
\end{equation}
which can be used to compute zeros in increasing order, starting at $z=z^{(0)}=0$. The first step
would be, because for $n$ odd $Y(0^+)/\dot{Y}(0^+)=0^+$ and for $n$ even 
$Y(0^+)/\dot{Y}(0^+)=-\infty$:
$$
z^{(1)}=T_{-1}(0)=\left\{
\begin{array}{l}
\Frac{2\pi}{\sqrt{L^2-4\lambda^2-1}},\quad n\, \mbox{odd},\\
\Frac{\pi}{\sqrt{L^2-4\lambda^2-1}},\quad n\, \mbox{even}.
\end{array}
\right.
$$

Once this first step is taken, we should compute $Y(z^{(1)})$ and $\dot{Y}(z^{(1)})$ and then proceed with the
next iteration. We propose the use of Taylor series for this computation. The difficulty 
in working in the $z$ variable 
is that Taylor series are not easy to implement: 
the successive derivatives of $Y(z)$ are not simple to compute by differentiation of the
ODE because the coefficients are no longer polynomials in $z$, and
therefore the derivatives don't satisfy a recurrence relation with a fixed number of terms. 
For this reason, we prefer to compute the functions 
by Taylor series in the 
original variable $x$, as considered in Section~\ref{taygen}. In addition, in order to avoid inversions 
of the variable in each step, we will write the fixed point method in the $x$ variable, although 
the underlying 
fixed point method will be (\ref{fixedz}) with $\Omega (z)$ given by (\ref{Yz}).

Then, working in the $x$ variable, we would start at:

\begin{equation}
\label{start}
x^{(1)}
=\tanh(T_{-1}(0))=\left\{
\begin{array}{l}
\tanh\left(\Frac{2\pi}{\sqrt{L^2-4\lambda^2-1}}\right),\, n\quad \mbox{odd},\\
\mbox{tanh}\left(\Frac{\pi}{\sqrt{L^2-4\lambda^2-1}}\right),\, n \quad \mbox{even}
\end{array}
\right.
\end{equation}
Then we compute $\tilde{Y}(x^{(1)})$ and $\tilde{Y}'(x^{(1)})$ using the 
Taylor series (\ref{truncs})  where, for this first step, 
$x=x^{(0)}=0$, $h=x^{(1)}-x^{(0)}=x_1$. We can take as
initial values
\begin{equation}
\label{tays}
\begin{array}{l}
\tilde{Y}(0)=\epsilon,\, \tilde{Y}'(0)=0,\quad n \mbox{ even, }\\
\tilde{Y}(0)=0, \tilde{Y}'(0)=\epsilon, \quad n \mbox{ odd, }
\end{array}
\end{equation}
with $\epsilon$ any fixed real number, say $\epsilon=1$
(as done in \cite{Gil:2019:FRA}, we can later renormalize the solutions using
one of the moments). 
 Once we have computed $\tilde{Y}(x^{(1)})$
and $\tilde{Y}'(x^{(1)})$ we iterate with the fixed point method, which we write in the $x$ variable using
$$
\Frac{\dot{Y}(z)}{Y(z)}=\Frac{dx}{dz}\Frac{dY}{dx}\Frac{1}{Y}=(1-x^2)\Frac{Y'(x)}{Y(x)},
$$
and in terms of the derivatives we are computing ($\tilde{Y}^{(k)}$) we have, because
$\tilde{Y}(x)=\sqrt{1-x^2}Y(x)$ (compare (\ref{yxz}) with (\ref{Ytil})),
\begin{equation}
\label{relse}
\Frac{Y(z)}{\dot{Y}(z)}=\Frac{\tilde{Y}(x)}{(1-x^2)\tilde{Y}'(x)+x \tilde{Y}(x)}.
\end{equation}
Given $x^{(1)}$ the next iteration of the fixed point method (\ref{fixedz}) is
\begin{equation}
\label{x1}
\mbox{arctanh}\,x^{(2)} =\mbox{arctanh}\,x^{(1)}-F(x^{(1)}),
\end{equation}
where
\begin{equation}
\label{Fx}
F(x)=\Frac{1}{\sqrt{\Omega (x)}}\arctan_{-1}\left(\sqrt{\Omega (x)}\Frac{\tilde{Y}(x)}{(1-x^2)\tilde{Y}'(x)+x\tilde{Y}(x)}\right)
\end{equation}
and 
\begin{equation}
\label{Ox}
\Omega(x)=\Frac{1}{4}\left[(L^2-1)(1-x^2)-4\lambda^2\right].
\end{equation}
We re-write (\ref{x1}) as
\begin{equation}
\label{ixtd}
x^{(2)}=g(x^{(1)})\equiv \Frac{x^{(1)}-\tanh(F(x^{(1)}))}{1-x^{(1)}\tanh(F(x^{(1)}))}.
\end{equation}

The algorithm proceeds similarly as described for Gauss--Hermite quadrature and, because $\Omega'(x)>0$ if
$x>0$, the nodes are computed in increasing order. Therefore, after the smallest positive
node $x_1$ has been computed by iterating $x^{(k+1)}=g(x^{k})$ (starting with (\ref{start}))
the next step would be
$$
x^{(1)}=g(x_1)=\Frac{x_1+\tanh\left(\pi/\sqrt{\Omega (x_1)}\right)}
{1+x_1 \tanh\left(\pi/\sqrt{\Omega (x_1)}\right)},
$$
and starting with this value we iterate again $x^{(k+1)}=g(x^{(k)})$ and compute the second positive node,
and so on. Parallel to this, the values of $\tilde{Y}(x^{(k+1)})$ and $\tilde{Y}'(x^{(k+1)})$ are computed
from Taylor series starting from $\tilde{Y}(x^{(k)})$ and $\tilde{Y}'(x^{(k)})$. At the same time the nodes
are computed, the values $\tilde{Y}'(x_i)$ are also obtained; from these values, we can compute the 
weights. 

We define the scaled weights as
$$
\omega_i =1/\tilde{Y}'(x_i)^2,
$$
and then the weights are given by (see (\ref{pesoz}))
\begin{equation}
\label{relw}
w_i =\gamma (1-x_i^2)^{\lambda} \omega_i,
\end{equation}
with $\gamma$ a factor which we can determine by normalizing to the moment of order $0$
$$
\Frac{\Gamma(\lambda +1)}{\Gamma (\lambda +3/2)}\sqrt{\pi}=\mu_0=
\displaystyle\int_{-1}^1 (1-x^2)^{\lambda} dx= 
\left[w_0\right]+2\displaystyle\sum_{i=1}^{\lfloor n/2 \rfloor}w_i,
$$
where the values $w_i$, $i>0$, are the weights corresponding to the 
positive nodes $x_1<x_2<\cdots x_{\lfloor n/2 \rfloor}$ and $w_0$ is the weight corresponding to the
node $x_0=0$ when $n$ is odd (which with the initial values (\ref{tays}) is $w_0=\gamma \omega_0 
=\gamma/\epsilon^2$);
the weight $w_0$ is inside brackets to denote that it only appears for $n$ odd.
Then, using (\ref{relw}),
\begin{equation}
\label{factor}
\gamma=\mu_0\left(\left[\Frac{\omega_0}{2}\right]+\displaystyle\sum_{i=1}^{\lfloor n/2 \rfloor}
(1-x_i^2)^{\lambda} \omega_i\right)^{-1},
\end{equation}
from which the weights (\ref{relw}) can be computed.

For analyzing the performance of this method, we have implemented the computation of 
Gauss-Gegenbauer quadratures both in Maple and Fortran (in double and quadruple precision arithmetics). 
In these programs, we normalize the weights to $2$ instead of $\mu_0$; the Gauss--Gegenbauer 
weights are then
recovered by multiplying them by the factor 
$\frac{\sqrt{\pi}}{2}\Gamma (\lambda+1)/\Gamma(\lambda+3/2)$, which is easily computed by Maple or 
with the Fortran program {\bf quotgamm} of the package {\bf gammaCHI} \cite{Gil:2015:GAP}.

As we later discuss, the algorithms are fast and accurate, except when $\lambda$ approaches $-1$. There
are two reasons for accuracy degradation in this case. In the first place, the largest zero tends to $1^{-}$
(and correspondingly the smallest zero to $-1^{+}$), which is problematic for Taylor series. 
In addition, 
as can be understood from the circle theorem, the most significant weight in this
case (and in general for $\lambda<-1/2$) is the last weight, corresponding to the largest node, and
 the errors in this last weight are carried to the rest of the weights due to
the final normalization step (\ref{factor}). In these cases, it is preferable not to compute the last weight
with Taylor series and to leave it instead as an unknown, to be fixed, together with the
normalization of the rest of weights, by using the first two even moments. With this, and denoting
$$
S_{0}=\left(\left[\Frac{\omega_0}{2}\right]+\displaystyle\sum_{i=1}^{\lfloor n/2 \rfloor -1}
(1-x_i^2)^{\lambda} \omega_i\right),\quad S_{x^2}=\displaystyle\sum_{i=1}^{\lfloor n/2 \rfloor -1}
x_i^2 (1-x_i^2)^{\lambda} \omega_i,
$$
we can compute 
\begin{equation}
\label{norwe}
\gamma=\Frac{x_{\lfloor n/2\rfloor}^2-(2\lambda+3)^{-1}}{x_{\lfloor n/2\rfloor}^2S_0-S_{x^2}},\quad 
w_{\lfloor n/2\rfloor}=1-\gamma S_{0}.
\end{equation}

We have implemented this additional step in our Fortran codes.
This correction reduces considerably the loss of accuracy as $\lambda\rightarrow -1$, 
as we will later discuss. 
However for the general case of
Jacobi quadrature we will discuss next, 
we prefer to recompute also the extreme nodes, and not only the weights, by using the
fixed point method associated to the angular change of variable (see Section~\ref{angular}).

\section{General Gauss--Jacobi quadrature}

The backbone of the general Gauss--Jacobi quadrature ($\alpha,\beta>-1$) will be again the fixed point method
 based on the transformation to
${\mathbb R}$, that is, the iteration (\ref{ixtd}), with $F(x)$ as in $(\ref{Fx})$ and with $\Omega(x)$ now given
by (\ref{Ozx}). The method of computation of $\tilde{Y}(x)$ and its derivative will be 
again the Taylor series of Section~\ref{erre} for the most part. There will be, however, exceptions to this.

In the first place, because in general $\alpha\neq \beta$, the problem is no longer symmetric around the
origin and we can not start with the initial values (\ref{tays}). The maximum of $\Omega (x)$ is 
placed at $x_e=(\beta^2 -\alpha^2)/(L^2-1)$, where  $L = 2n + \alpha + \beta + 1$,  and we should start
at this point, computing $\tilde{Y}(x_e)/\tilde{Y}'(x_e)$ in order to start the process. We will compute 
this starting value with the three-term recurrence relation (\ref{TTRR}), as we explain in the next subsection. In the second
place, as already described for the Gegenbauer case, the extreme nodes need particular attention, and we will
recompute them using the angular transformation of Section~\ref{angular}, with the functions computed via the 
continued fraction (\ref{cfP}).

In our algorithms we use the fact that $P^{(\alpha,\beta)}(x)=P^{(\beta,\alpha)}(-x)$ and so, instead of performing
a forward sweep for $x>x_e$ and a backward sweep for $x<x_e$, we perform two forward sweeps: one for the original 
values of $\alpha$ and $\beta$ and a second one with interchanged values (the signs of the nodes are changed after this
second computation).  

\subsection{Recurrence relation}

As mentioned, we start the process computing $\tilde{Y}(x_e)/\tilde{Y}'(x_e)$ by using the three-term
recurrence relation. For this purpose, and in order to avoid overflows, it is better to re-write the 
recurrence relation (\ref{TTRR}) in terms of ratios as follows

\begin{equation}
\label{recrat}
\begin{array}{l}
\Frac{P_{n+1}^{(\alpha,\beta)}(x)}{P_{n}^{(\alpha,\beta)}(x)}=
\Frac{1}{A_n}\left[B_n -\Frac{C_n}{P_{n}^{(\alpha,\beta)}(x)/P_{n-1}^{(\alpha,\beta)}(x)}\right],\\
\\
A_n=2(n+1)(n+\alpha+\beta+1)(L-1),\\
B_n=L\left\{(L^2-1)x+\alpha^2-\beta^2\right\},\\
C_n=2(L+1)(n+\alpha)(n+\beta),
\end{array}
\end{equation}
with starting value $P_{1}^{(\alpha,\beta)}(x)/P_{0}^{(\alpha,\beta)}(x)=(\alpha-\beta+
(\alpha+\beta+2)x)/2$.

With the notation used so far, $\tilde{Y}(x)=(1-x)^{(\alpha +1)/2} (1+x)^{(\beta +1)/2} P_n^{(\alpha,\beta)}(x)$ 
and considering the
derivative rule \cite[18.9.17]{Koorn:2010:OP} we have

\begin{equation}
\label{yrat}
\begin{array}{lll}
\Frac{\tilde{Y}'(x)}{\tilde{Y}(x)}&=&\Frac{\beta +1}{2(1+x)}-\Frac{\alpha +1}{2(1-x)}+
\Frac{P_{n}^{(\alpha,\beta)\prime}(x)}{P_{n}^{(\alpha,\beta)}(x)}\\
&=&\Frac{n+\beta +1}{2(1+x)}-\Frac{n+\alpha +1}{2(1-x)}\\&&+
\Frac{1}{(L-1)(1-x^2)}\left\{n(\alpha-\beta)+2(n+\alpha)(n+\beta)\Frac{P_{n-1}^{(\alpha,\beta)}(x)}{
P_{n}^{(\alpha,\beta)}(x)}\right\}.
\end{array}
\end{equation}

Then, combining (\ref{recrat}) and (\ref{yrat}), the ratio $\tilde{Y}(x_e)/\tilde{Y}'(x_e)$ can be computed in
order to start the process. The next step would be to compute the first iteration with 
(\ref{ixtd}),  with $\Omega(x)$ given
by (\ref{Ozx}).

\subsection{Extreme nodes: angular variable and continued fraction}
\label{extreme}

For the extreme nodes, and particularly for computing the weights, the angular variable $x=\cos\theta$ 
is more convenient. Let us, for instance, consider the computation of the zeros close to $x=1$. Assume 
that one of such zeros is $x_i=1-\delta$, with $\delta$ a small number, then 
$\delta=1-x_i=2\sin^2 (\theta_i/2)$,
and if the value $\theta_i$ is determined with a given relative precision, then the relative 
accuracy of $x=\cos\theta_i=1-2\sin^2 (\theta_i/2)$ will be higher. 
On the other hand, the attainable accuracy of the corresponding 
weights will be also higher by using (\ref{defh}) in the angular variable instead of Taylor series
in the original variable, which are problematically 
close to the singularities of the ODE.

We will use the angular change of variables only to refine the extreme nodes and weights already
computed with the change $x=\tanh z$. We do this for at least the three largest zeros when $\alpha<0$
and similarly for the smallest negative nodes when $\beta<0$; as $n$ increases, we increase logarithmically
the number of extreme zeros computed in this alternative way (one more zero as the degree increases by one 
order of magnitude).  These negative parameter cases are indeed 
the most problematic ones because the largest node tends to $+1$ when $\alpha\rightarrow -1^{+}$ and the smallest 
node to $-1$ as $\beta\rightarrow -1^{+}$. For non-negative parameters, 
the methods described in this section are not needed. 

Considering the angular transformation $x=\cos\theta$ of Section~\ref{angular}, 
we have that $Y(\theta)$, given by Eq.~(\ref{yxa}), satisfies the second order
ODE $\ddot{Y}(\theta)+\Omega (\theta) Y(\theta)=0$, with $\Omega(\theta (x))$ given by (\ref{eyx}). In the
$\theta$ variable this reads
$$
\Omega (\theta)=\Frac{1}{4}\left[L^2+\Frac{\frac14-\alpha^2}{\sin^2 (\theta/2)}
+\Frac{\frac14-\beta^2}{\cos^2 (\theta/2)}\right].
$$

The starting values for the fixed point method are the estimations given by the principal method
 (based on the change $x=\tanh z$), and the new iteration is used to improve such values. We then
start from $\theta_i =\arccos (x_i)$ where $x_i$ are the extreme zeros computed with the
 principal method. In this case it is convenient
to use the fixed point iteration (\ref{fixedpoint2}) with bilateral convergence, that is

$$
g(\theta)=\theta-\Frac{1}{\sqrt{\Omega (\theta)}}\arctan\left(\sqrt{\Omega (\theta)}Y(\theta)/\dot{Y}
(\theta)\right).
$$

For computing the ratio $Y(\theta)/\dot{Y}(\theta)$ we use the definition of $Y(\theta)$ 
together with (\ref{dera}) and (\ref{reca}); we have

$$
\sin\theta\Frac{\dot{Y}(\theta)}{Y(\theta)}\equiv h(\theta)^{-1}=1/2+\alpha+L\sin^2 \Frac{\theta}{2}
-2(n+\alpha+\beta+1)\sin^2\Frac{\theta}{2}\Frac{P_n^{(\alpha+1,\beta)}(x)}{P_n^{(\alpha,\beta)}(x)},
$$
where the ratio of Jacobi polynomials can be computed with the continued fraction 
(\ref{cfP}), conveniently written in the variable $\theta$ (replacing $1-x$ by $2\sin^2(\theta/2)$).

Then the fixed point method can be written

$$
g(\theta)=\theta -\Frac{\sin\theta}{\sqrt{\Delta}}\arctan \left(\sqrt{\Delta} h(\theta) \right)
$$
with
$$
\Delta=\Frac{1}{4}-\alpha^2+(\alpha^2-\beta^2)
\sin^2 \Frac{\theta}{2}+\Frac{L^2}{4}\sin^2\theta .
$$

Once the extreme nodes have been refined in the angular variable, the weights can be refined too. For this
purpose, we consider the expression of the weights and the relation of Jacobi polynomials with hypergeometric
functions, which, together with the derivative rule for Gauss hypergeometric functions leads to

\begin{equation}
\label{extewe}
w_i= \Frac{K_{n,\alpha,\beta}}{\sin^2\theta_i\,
\left[{}_2 {\rm F}_{1} \left(-n+1, n+\alpha+\beta+2;\alpha+2;
\sin^2\left(\Frac{\theta_i}{2}\right)\right)\right]^2},
\end{equation}
where $K_{n,\alpha,\beta}$ is a constant not depending on $\theta_i$ which can be obtained from 
Eqs. (\ref{peso}) and (\ref{defh})
$$
K_{n,\alpha,\beta}=\left(\Frac{2(n-1)!}{(n+\alpha+\beta+1)(\alpha+2)_{n-1} }\right)^2 M_{n,\alpha,\beta},
$$
with $M_{n,\alpha,\beta}$ given by (\ref{constante}). 

Because the argument of the terminating series $\sin^2\left(\Frac{\theta_i}{2}\right)$ will be small
for the extreme zeros, few terms of this series will be needed for an 
accurate computation close to $x=1$, also for large $n$.

It is possible to skip the computation of these
constants, in the same way that for the symmetric case $\alpha=\beta=\lambda$ we didn't need to compute 
$M_{n,\lambda,\lambda}$.  
For fixing the normalization of the weights we should take into account that we may have
up to three sets of weights with different normalizations: the weights computed by Taylor series (principal method) 
and up to
two sets of extreme zeros (the positive and the negative), which are computed independently. 
One possibility to fix the normalizations is to use the first three moments of the weight, that is
using that the $n$ nodes and weights satisfy (for $n\ge 2$):
\begin{equation}
\begin{array}{l}
\displaystyle\sum_{i=1}^n w_i =\mu_0=\displaystyle\int_{-1}^{1}(1-x)^{\alpha}(1+x)^{\beta}dx = 
2^{\alpha+\beta+1}\Frac{\Gamma (\alpha+1)\Gamma (\beta+1)}{\Gamma (\alpha+\beta+2)},\\
\displaystyle\sum_{i=1}^n x_i w_i =\mu_1=\displaystyle\int_{-1}^{1}x(1-x)^{\alpha}(1+x)^{\beta}dx 
=\mu_0\Frac{\beta-\alpha}{\alpha+\beta+2},\\
\displaystyle\sum_{i=1}^n x_i^2 w_i =\mu_2=\displaystyle\int_{-1}^{1}x^2(1-x)^{\alpha}(1+x)^{\beta}dx 
=\mu_0\Frac{(\alpha-\beta)^2+\alpha+\beta+2}{(\alpha+\beta+2)(\alpha+\beta+3)}.
\end{array}
\end{equation}

We have observed, however, that the resulting linear system for these normalization constants loses some accuracy when either
$\alpha$, $\beta$ or both are close to $-1/2$.  It is,
however, very accurate for parameters close to $-1$, when the extreme zeros are the dominant ones.

As an alternative to avoid inaccuracies for parameters close to $-1/2$ 
\footnote{It is interesting to observe that, as we discuss later, also the Golub-Welsch algorithm appears to suffer
from some loss of accuracy for these parameter values}
we compute the extreme weights with 
formula (\ref{extewe}), with the constant computed in terms of gamma functions\footnote{Observe that we need
 to compute $K_{n,\alpha,\beta}$ only if $\alpha<0$ (and $K_{n,\beta,\alpha}$ only if $\beta<0$), and in this
case, in practical terms the constant does not overflow/underflow as $n$ becomes large (it does so algebraically). 
However, the gamma functions do become
huge. For computing this it is preferable to compute the logarithm of the constant and exponentiate afterwards. The
 logarithm of the gamma function is a widely available computation, for example through the command {\bf gammaln}
in Matlab or with the Fortran function {\bf loggam} of \cite{Gil:2015:GAP}}. The only normalization to be 
determined is
for the weights computed with Taylor series, which we can determine with the moment of order zero. 
Then, if, say, $\{w_i\}_{i=1}^{n_l}$ and $\{w_i\}_{i=n_u}^{n}$ are the nodes computed in the angular variable
(which are final weights, with no normalization required), and $\{\tilde{w}_i\}_{i=n_{l}+1}^{n_{u}-1}$ are
the weights computed by Taylor series, related with the final weights by $w_i=\gamma \tilde{w}_i$ we 
have
$$
\mu_0 =\gamma \tilde{S}_{T} + S_{\theta},\,\tilde{S}_T =\displaystyle\sum_{i=n_l+1}^{n_u-1}\tilde{w}_i,\,
S_{\theta}=\displaystyle\sum_{i=1}^{n_l}w_i + 
\displaystyle\sum_{i=n_u}^{n}w_i,
$$
from where we can compute $\gamma$, and then all the final weights $\{w_i\}_{i=1}^n$ 
are obtained.

In our Matlab codes we adopt this scheme when both $\alpha$ and $\beta$ are larger than $-3/4$ and switch to 
the approach in terms of the
three first moments in the other case.

This ends the description of the methods used for the computation of Gauss--Jacobi quadrature. 
Next we describe the performance of the resulting algorithms.

\section{Numerical tests}

We now test the several implementations of our algorithms. We start by describing the 
high-accuracy performance of our methods, in particular for Gauss--Gegenbauer quadrature, 
and later we compare, both in terms of speed and accuracy, our double precision version of 
our algorithms for the general Jacobi case against the Chebfun \cite{Dris:2014:CG} implementation of the 
methods described in 
\cite{Bogaert:2014:IFC,Hale:2013:FAA} and against the
Golub--Welsch algorithm \cite{Golub:1969:COG}. Finally, we perform some additional tests for the 
general Jacobi case by comparing our algorithms
with a high accuracy implementation of the Golub-Welsch algorithm using Maple.

\subsection{Gauss--Gegenbauer quadrature and high-accuracy computations}

The fact that the methods are based on convergent processes and that the nonlinear method is of order 
four makes this an ideal approach for arbitrary accuracy computations. With this method,
it is possible to efficiently compute high order quadrature rules with high accuracy. Furthermore, 
as the degree is higher the cost of computation per node becomes smaller, both in terms of the number of
iterations per node and the number of Taylor sums per node. 

In order to illustrate
these facts, we have implemented in Maple the core method (based solely on the $x=\tanh z$ 
transformation) for Gauss--Gegenbauer quadratures and for increasing
degrees and accuracies. These results are illustrated in Table \ref{tab1}.

\begin{table}
$$
\begin{tabular}{c|ccccccc}
n \textbackslash D & 16 & 32 & 64 & 128 & 256 & 512 & 1024 \\
\hline
 10    & 3/367 & 5/556 & 4.4/787 & 5/954 & 5.4/1203 & 6/1360 & 6.4/1582\\ 
 100   & 2.8/118 & 3.2/158 & 4/222 & 4.3/325 & 5/510 & 5.4/778 & 6/1009\\ 
 1000   & 2/78 & 3/114 & 3.5/156 & 4/229 & 4.7/359 &  5/581 & 5.7/819 \\
10000   & 2/75 & 3/110 & 3/139 & 4/214 & 4.1/319 & 5/532  & 5.1/732 \\
100000 &  2/75 & 2.5/100 & 3/137 & 3.9/207 & 4/307 & 5/507 & 5/686 \\
\end{tabular}
$$
\caption{Average number of iterations per node and number of the terms of the Taylor series used 
per node for degrees 
$n=10,\,100,\,1000,\,10000,\,100000$ (rows) and for relative accuracies of $10^{-D}$ 
(columns). In this table $\lambda=-0.8$.}%
\label{tab1}%
\end{table}

As we can observe the number of iterations per node, even for extreme accuracies,
 usually does not exceed $6$, of course increasing as the accuracy increases, and with a rate 
corresponding to a fourth order fixed point method (roughly one more iteration when the number of digits
is quadrupled). This behavior is observed regardless of the value of $\lambda=\alpha=\beta$.

\begin{figure}[tb]
\vspace*{0.8cm}
\begin{center}
\begin{minipage}{3cm}
\centerline{\includegraphics[height=4cm,width=6cm]{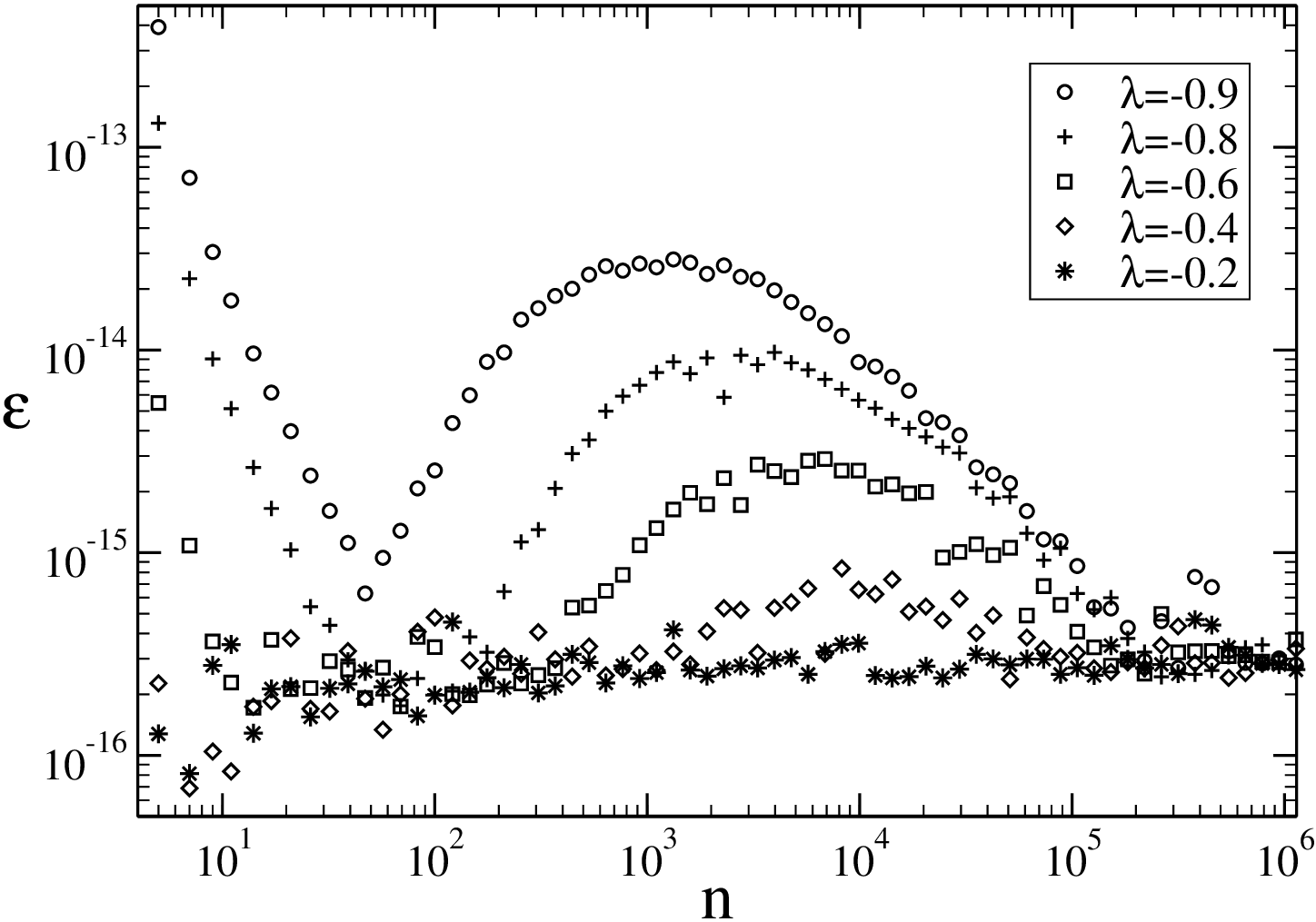}}
\end{minipage}
\hspace*{3cm}
\begin{minipage}{3cm}
\centerline{\includegraphics[height=4cm,width=6cm]{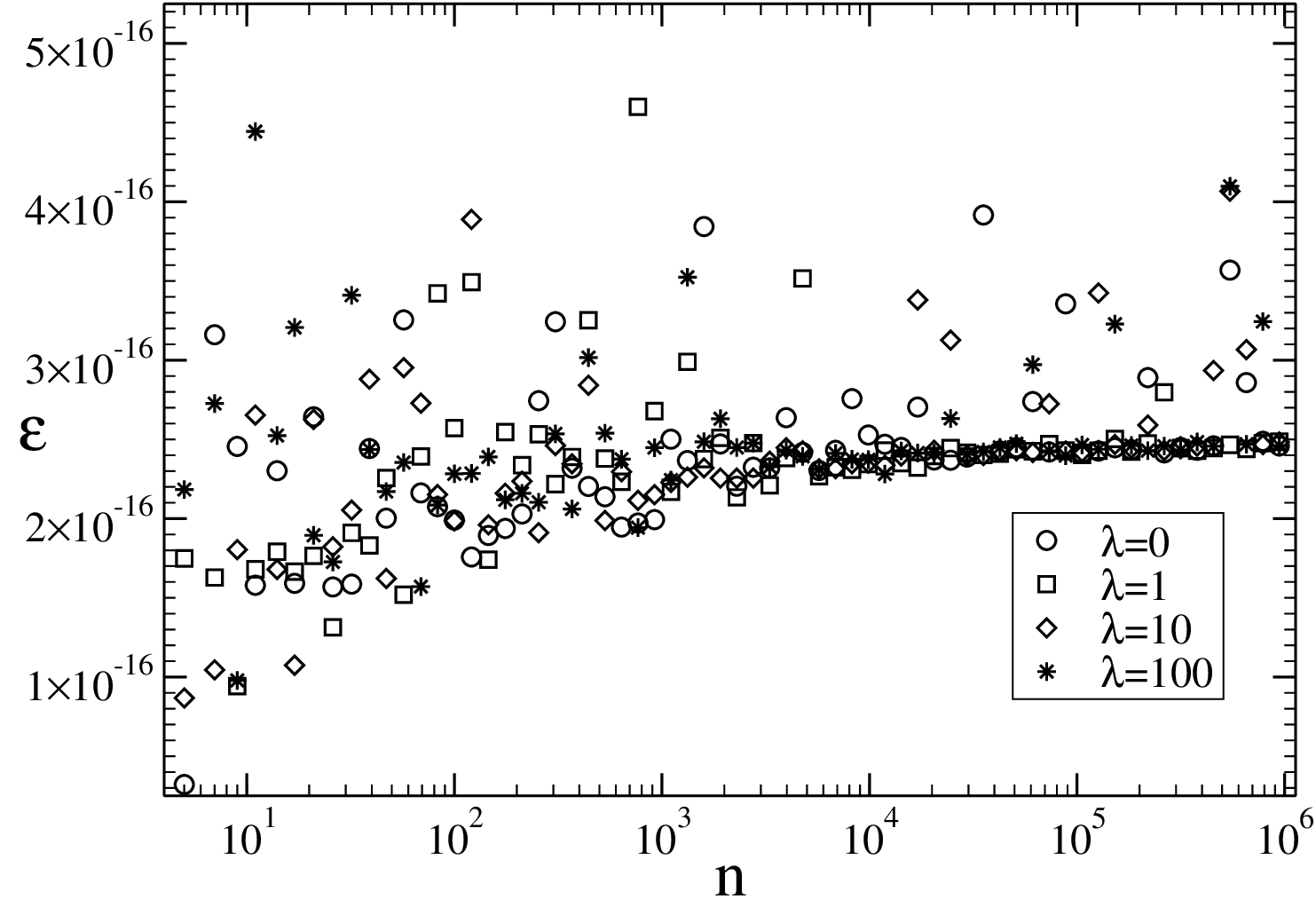}}
\end{minipage}
\end{center}
\caption{Maximum relative error in the computation of the Gauss--Gegenbauer nodes for negative values
of $\lambda=\alpha=\beta$ (left) and non-negative values (right).}
\label{fig1}
\end{figure}

\begin{figure}[tb]
\vspace*{0.8cm}
\begin{center}
\begin{minipage}{3cm}
\centerline{\includegraphics[height=4cm,width=6cm]{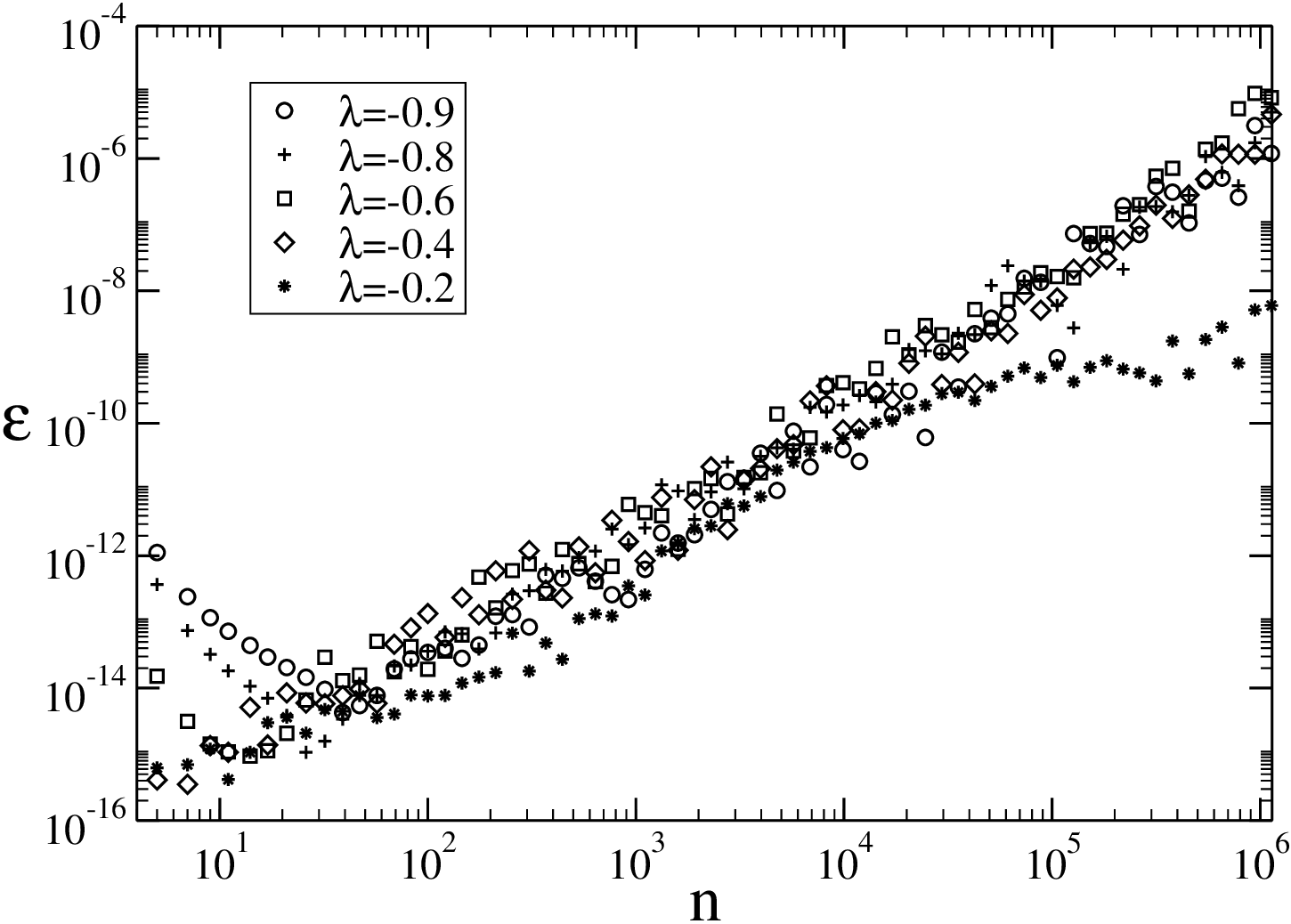}}
\end{minipage}
\hspace*{3.2cm}
\begin{minipage}{3cm}
\centerline{\includegraphics[height=4cm,width=6cm]{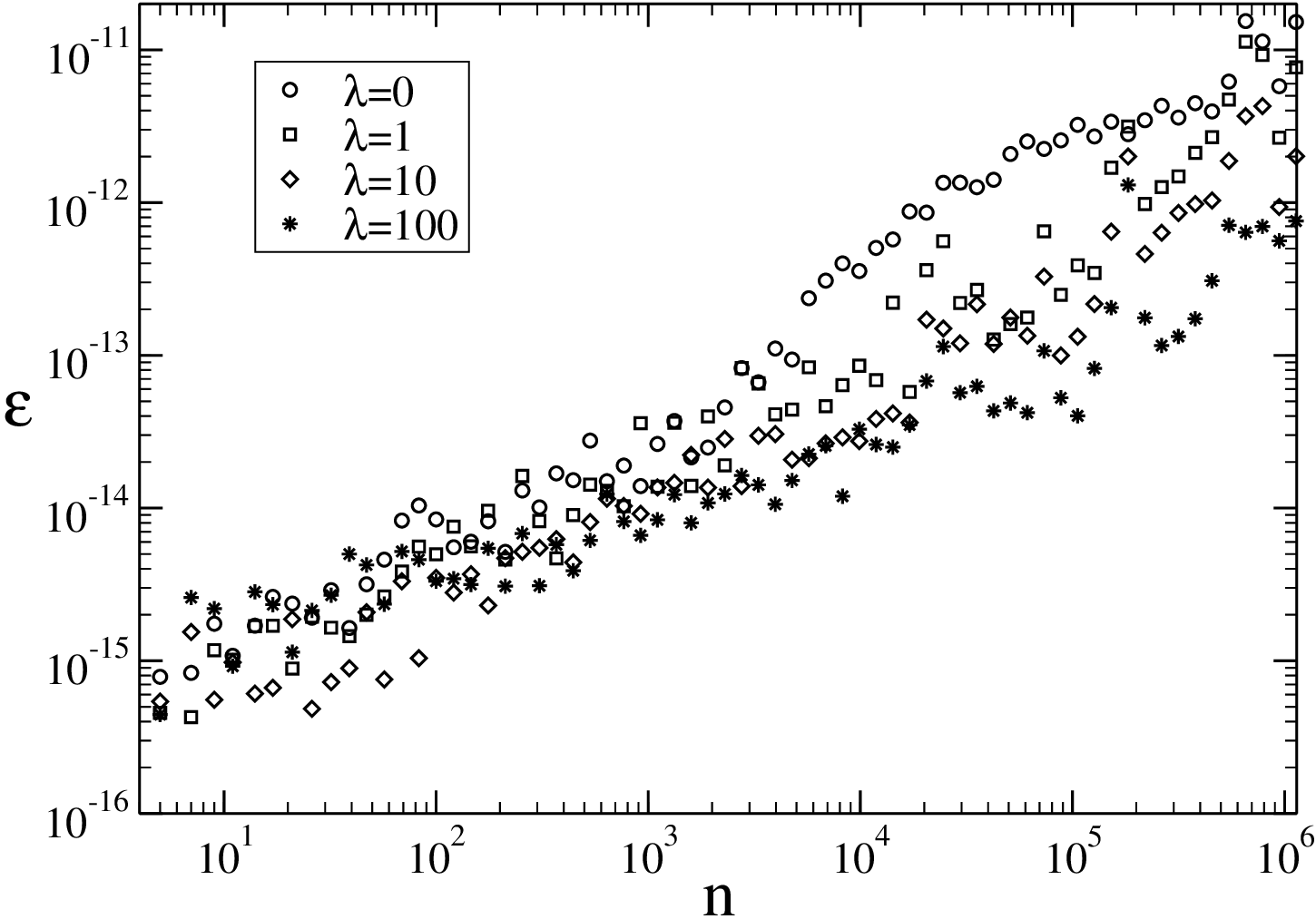}}
\end{minipage}
\end{center}
\caption{Relative maximum error $\varepsilon_{rm}(\{w\})$ in the computation of the Gauss--Gegenbauer weights for negative values
of $\lambda=\alpha=\beta$ (left) and non-negative values (right). The double precision weights $w_i^{d}$ 
are compared against the quadruple precision ones $w_i^{q}$.}
\label{fig2}
\end{figure}

For numerically testing the accuracy, we have checked the consistency of the computation of the nodes and
weights with different accuracies. In order to facilitate these tests and being able 
to perform more intensive ones, we have translated the Maple algorithm to Fortran 90, both in double and 
quadruple precision \footnote{The Maple worksheet (Gauss-Gegenbauer) and 
the Matlab code (Gauss-Jacobi) mentioned in this paper can be found at 
https://personales.unican.es/segurajj/gaussian.html, together
with the codes corresponding to the Gauss--Hermite and Gauss--Laguerre cases of \cite{Gil:2019:FRA}. None 
of these codes should be considered as final
versions of our algorithms.}. In Figures~\ref{fig1} and \ref{fig2}, 
we compare the output of the
double precision implementation against the quadruple versions.

In Fig. \ref{fig1} we plot the maximum relative error of the positive nodes, comparing the nodes in double
precision $x_i^{(d)}$ with the same nodes in quadruple precision $x_i^{(q)}$. That is, we plot, as a function
of $n$ and for different values of $\lambda$, the quantity
\begin{equation}
\varepsilon_{mr}(\{x\})=\displaystyle\max_{i} \varepsilon_r (x_i)= 
\displaystyle\max_{i}\left|1-\Frac{x_i^{(d)}}{x_i^{(q)}}\right|,
\end{equation}
where $\varepsilon_r (x_i)$ is the relative error for the node $x_i$.

In Fig. 1 left, the maximum errors for the nodes are shown for negative values of $\lambda$. We observe that the maximum relative errors
are close to double precision accuracy except when $\lambda <-0.5$; this is due, as commented in Section~\ref{GGQ},
 to the fact that Taylor series lose some precision for the extreme nodes as 
$\lambda \rightarrow -1^{+}$,
and the largest relative errors take place for the extreme zeros. As we discuss later, this loss of 
accuracy for 
the nodes is solved by considering the angular variable for few of the extreme nodes, as described in 
Section~\ref{extreme}.

 Fig. \ref{fig2} shows the
relative maximum (absolute) error for the weights, that is
\begin{equation}
\varepsilon_{rm}(\{w\})=
\displaystyle\max_{i}\varepsilon_{ra}(w_i)=
\Frac{\displaystyle\max_{i}|w_i^{(d)}-w_i^{(q)}|}{\displaystyle\max_{i}w_i^{(q)}},
\end{equation}
where $\varepsilon_{ra}(w_i)=|w_i^{(d)}-w_i^{(q)}|/(\max_{i}w_i^{(q)})$ 
is the absolute error for the weight $w_i$ relative to the maximum weight 
%\footnote{Absolute errors  do not make much sense as an error measure for the weights, 
%because the magnitude of the weights depends strongly on the parameters and the 
%corresponding nodes, particularly for large values of the parameters (see Eq. (\ref{circle})). 
%It is more reasonable to consider absolute errors relative to the maximum value of the weight, as we do here.}.
Relative maximum error was also used as error measure in this same context in \cite{Glaser:2007:AFA,Hale:2013:FAA}. 
 This is a reasonable measure for the weights because in the evaluation of quadrature rules, 
the largest weights are the most significant ones, while if a weight is much smaller
than the maximum it should be enough to compute it with lower relative accuracy; some
weights may be even smaller than the underflow number (for very large parameters $\alpha$ and or $\beta$). 
We also consider relative error for the weights later when we test the algorithm for the general Jacobi case. 

In Fig. \ref{fig2} we observe that there is a gradual loss of accuracy as the degree increases. This loss of accuracy was also observed
in \cite{Glaser:2007:AFA} for Gauss--Legendre, and it is surely related to the successive application of Taylor series, which was
also used in \cite{Glaser:2007:AFA}. Still, for non-negative $\lambda$ the relative maximum error is close to $10^{-14}$
even for degrees as large as $1000$. Same as happened with the nodes, there is some additional loss of accuracy for the
 weights corresponding to 
the extreme nodes
when $\lambda$ is negative, which is in part corrected by using (\ref{norwe}). As we will see for the more general algorithm 
(that we have implemented in Matlab), the use of the angular variable as described in Section~\ref{extreme} will improve the 
accuracy.

We recall that asymptotic approximations are accurate for degrees $n\ge 100$ and
$-1<\lambda\le 5$ \cite{Gil:2019:NIC}, and with close to double relative accuracy for both the nodes and the weights, 
and for that cases
such approximations are preferable for double precision computations. However, outside this range or when higher accuracy
is needed, the algorithm presented in this paper is the best option.

\subsection{General Gauss--Jacobi algorithm}

We have implemented our algorithm for the general Jacobi case in Matlab, and we plan to implement this
in Fortran and Maple in the near future. One of the reasons to choose this platform is that there are
some alternative methods to compare with, in particular the chebfun \cite{Dris:2014:CG} 
implementation of the methods in \cite{Hale:2013:FAA,Bogaert:2014:IFC} 
and the classical Golub--Welsch algorithm \cite{Golub:1969:COG}, which can be easily programmed 
in Matlab using its powerful matrix diagonalization
routines. We will compare our method (which we label as NEW) 
against the chebfun program {\bf jacpts.m} for computing Gauss--Jacobi quadrature (labeled as CHEB), and our own implementation of the Golub--Welsch algorithm (GW). These three methods, implemented in double precision accuracy in Matlab, 
allows us to perform quite extensive tests as a function 
of the degree $n$ and/or the parameters, particularly when comparing CHEB with NEW, which are quite efficient methods; the comparison
 with GW is more
time consuming and sets a limit on the value of $n$. The conclusions that will be drawn from these tests will
 be also corroborated for some specific values of $n$, $\alpha$ and $\beta$ by comparing the results with a high
 accuracy computation of the nodes and weights using a variable precision implementation of the GW algorithm (in Maple).

\subsubsection{The symmetric case}

\begin{figure}[tb]
\vspace*{0.8cm}
\begin{center}
\begin{minipage}{3cm}
\centerline{\includegraphics[height=4cm,width=6cm]{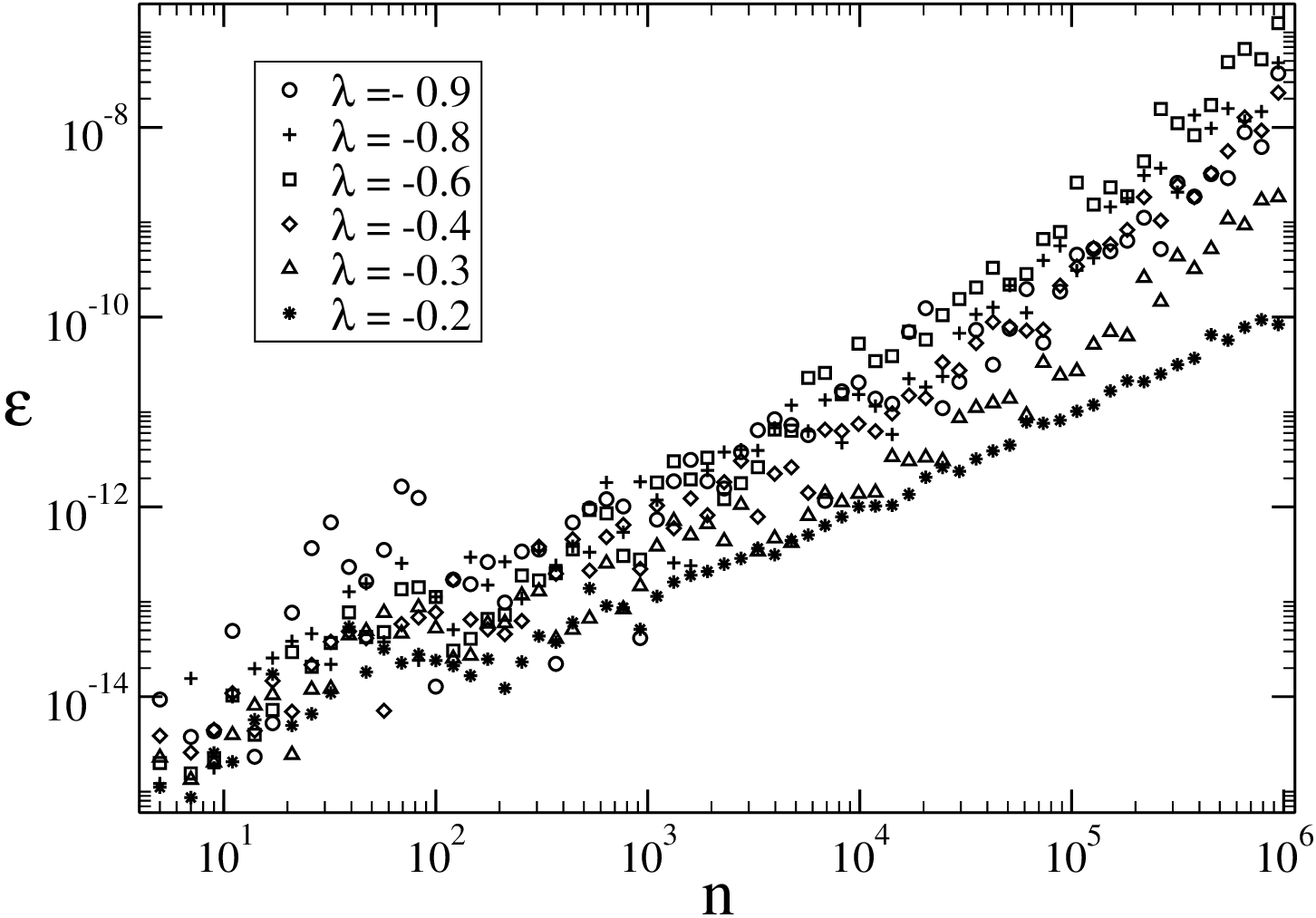}}
\end{minipage}
\hspace*{3.2cm}
\begin{minipage}{3cm}
\centerline{\includegraphics[height=4cm,width=6cm]{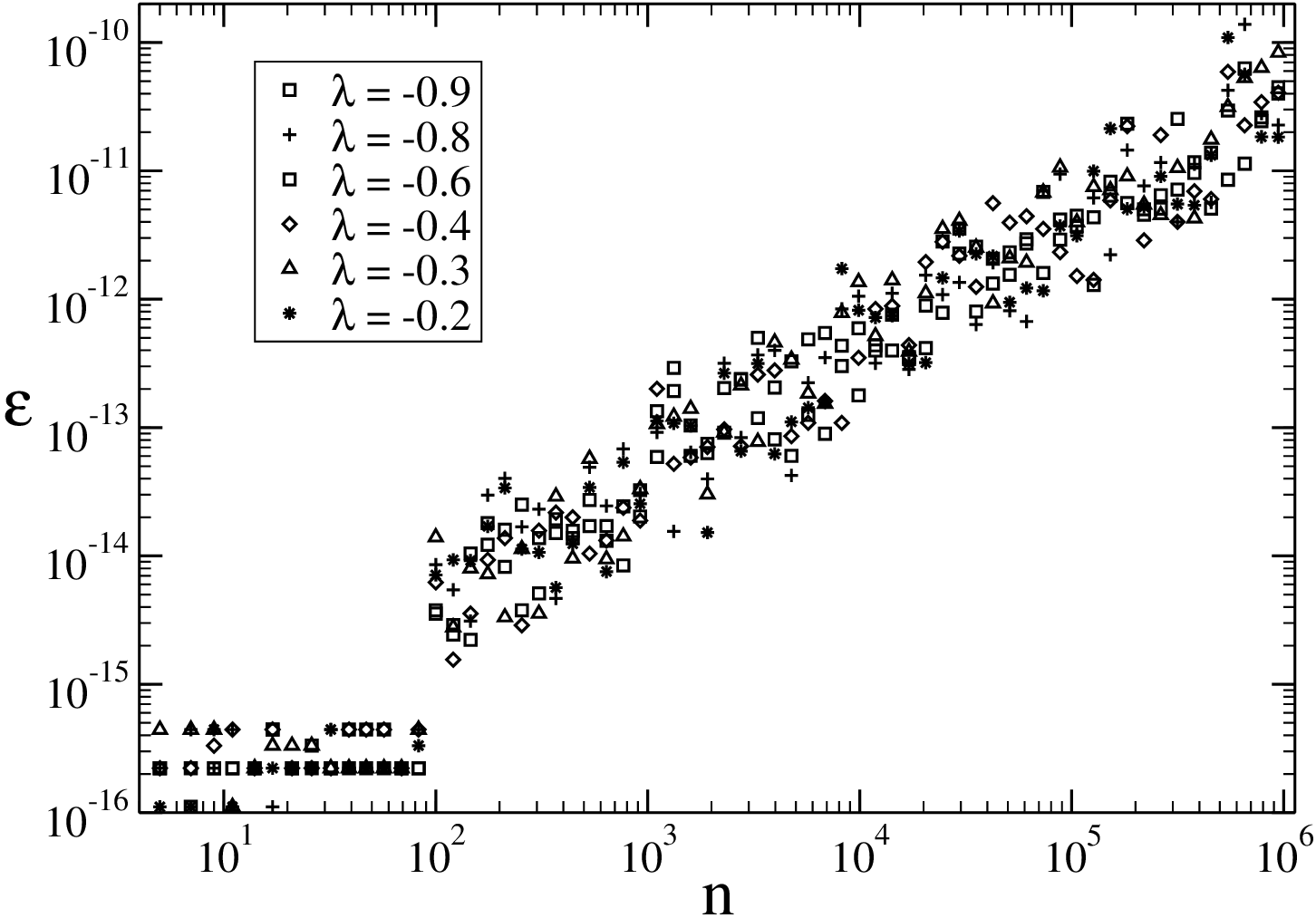}}
\end{minipage}
\end{center}
\caption{Error in the computation of the Gauss--Gegenbauer nodes and weights of our Matlab implementations against
the chebfun algorithm {\bf jacpts}, implementing the methods in \cite{Hale:2013:FAA}. Left: relative maximum error for
the weights. Right: maximum relative error for the nodes.}
\label{fig3}
\end{figure}

We start our comparison by first restricting to the symmetric case $\alpha=\beta=\lambda$.
In Fig. \ref{fig3} we show the maximum relative errors for the nodes and the relative maximum error for the weights
as a function of $n$ and for various negative values of $\lambda$. 
Non-negative values are not considered because
for that case the computation is exactly as in the previously discussed Fortran implementation, and the errors for our algorithm 
are those shown in Fig. \ref{fig1} and Fig. \ref{fig2} (right).

In Fig. \ref{fig3} we notice a difference in the results for $n<100$ and $n\ge 100$
due to the fact that the algorithm {\bf jacpts} uses different methods on those two cases: the polynomials are evaluated by
the three-term recurrence relation for $n<100$ and with asymptotics otherwise.
As we can observe, the accuracy worsens as $n$ increases both for the nodes and the weights (except for the nodes when $n<100$). 
For
the case of the weights, and as before discussed, this error degradation is due to the new method, while for the nodes 
it is due to CHEB, because it computes the nodes in the angular variable, and when inverting to compute the nodes $x_i$, relative accuracy
is not kept for the nodes close to zero. For the weights when $n<100$ we also observe, particularly for the case $\lambda=-0.9$,
 some errors larger than the rest. By comparing with the GW algorithm the largest error corresponds to the  
CHEB algorithm. When $n<100$ the nodes are correct within double precision accuracy for both CHEB and NEW. The main conclusion,
apart from particular behaviors for $n<100$, is that there is error degradation as $n$ increases for the nodes for CHEB,
 for the weights for NEW and for both the nodes and the weights for GW (not shown). GW, in addition, becomes prohibitively slow for
large $n$.

\subsubsection{The general case: comparing three methods in double precision}

\begin{figure}[tb]
\vspace*{0.8cm}
\begin{center}
\begin{minipage}{3cm}
\centerline{\includegraphics[height=4cm,width=6cm]{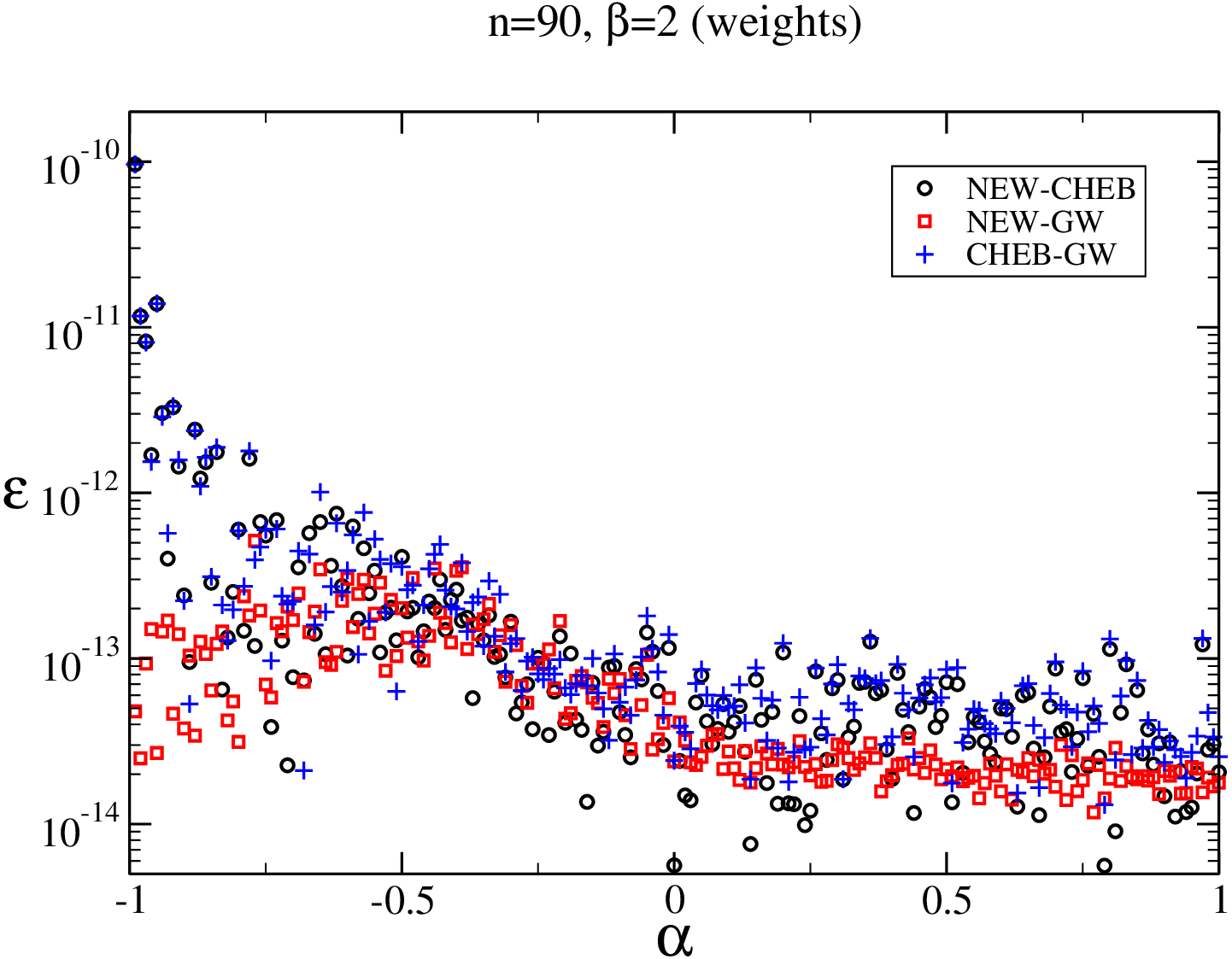}}
\end{minipage}
\hspace*{3.2cm}
\begin{minipage}{3cm}
\centerline{\includegraphics[height=4cm,width=6cm]{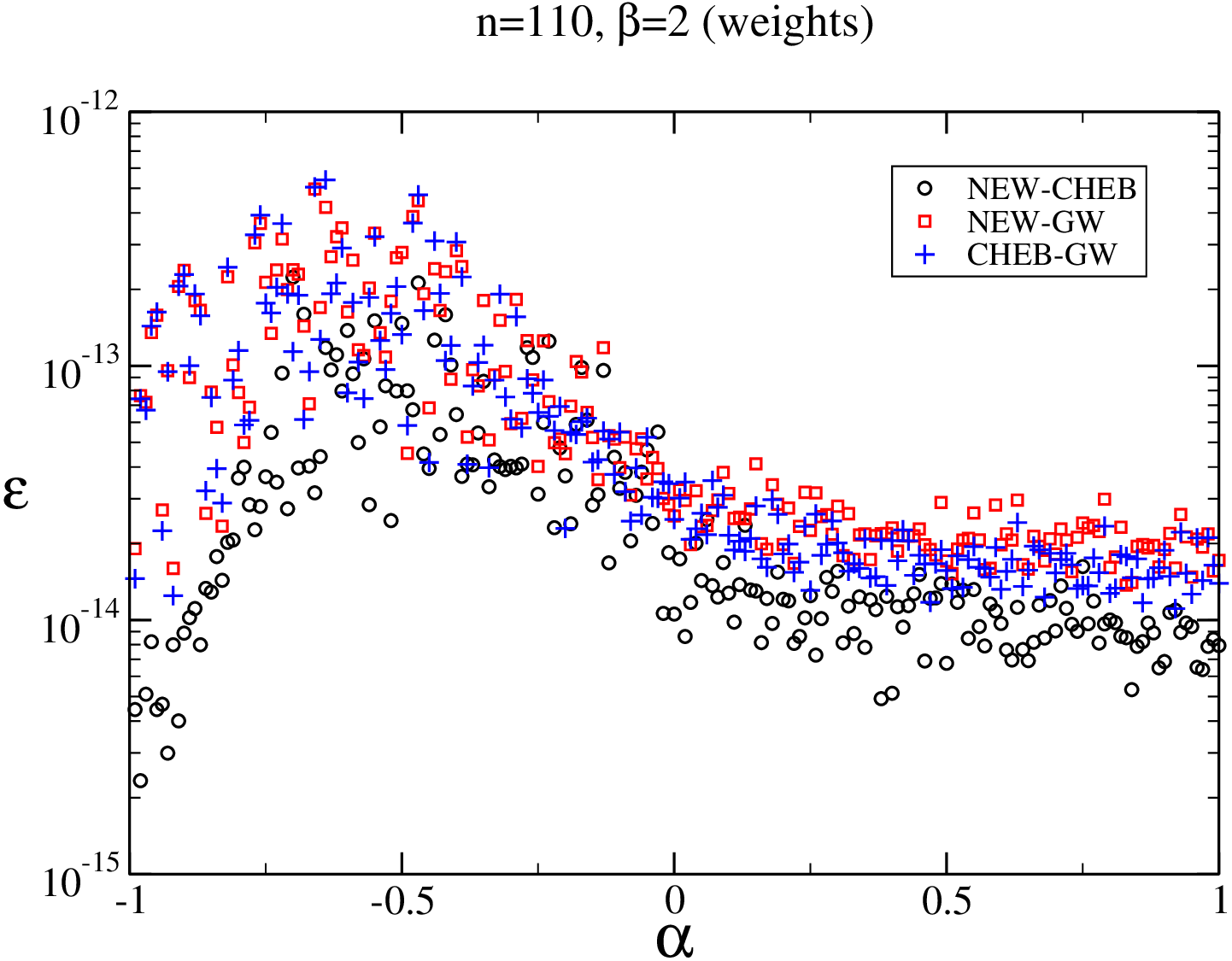}}
\end{minipage}
\end{center}
\vspace*{0.2cm}
\begin{center}
\begin{center}
\begin{minipage}{3cm}
\centerline{\includegraphics[height=4cm,width=6cm]{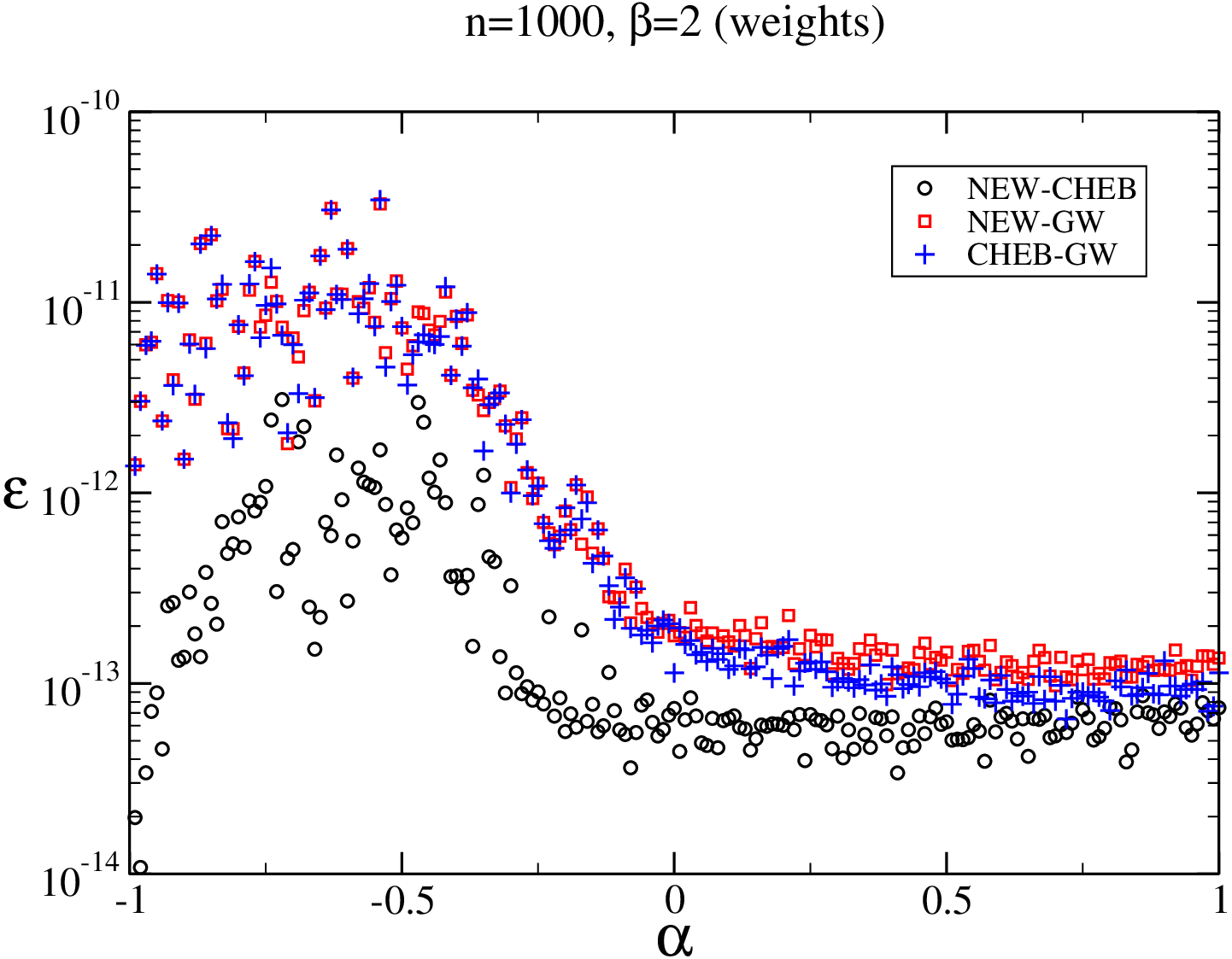}}
\end{minipage}
\hspace*{3.2cm}
\begin{minipage}{3cm}
\centerline{\includegraphics[height=4cm,width=6cm]{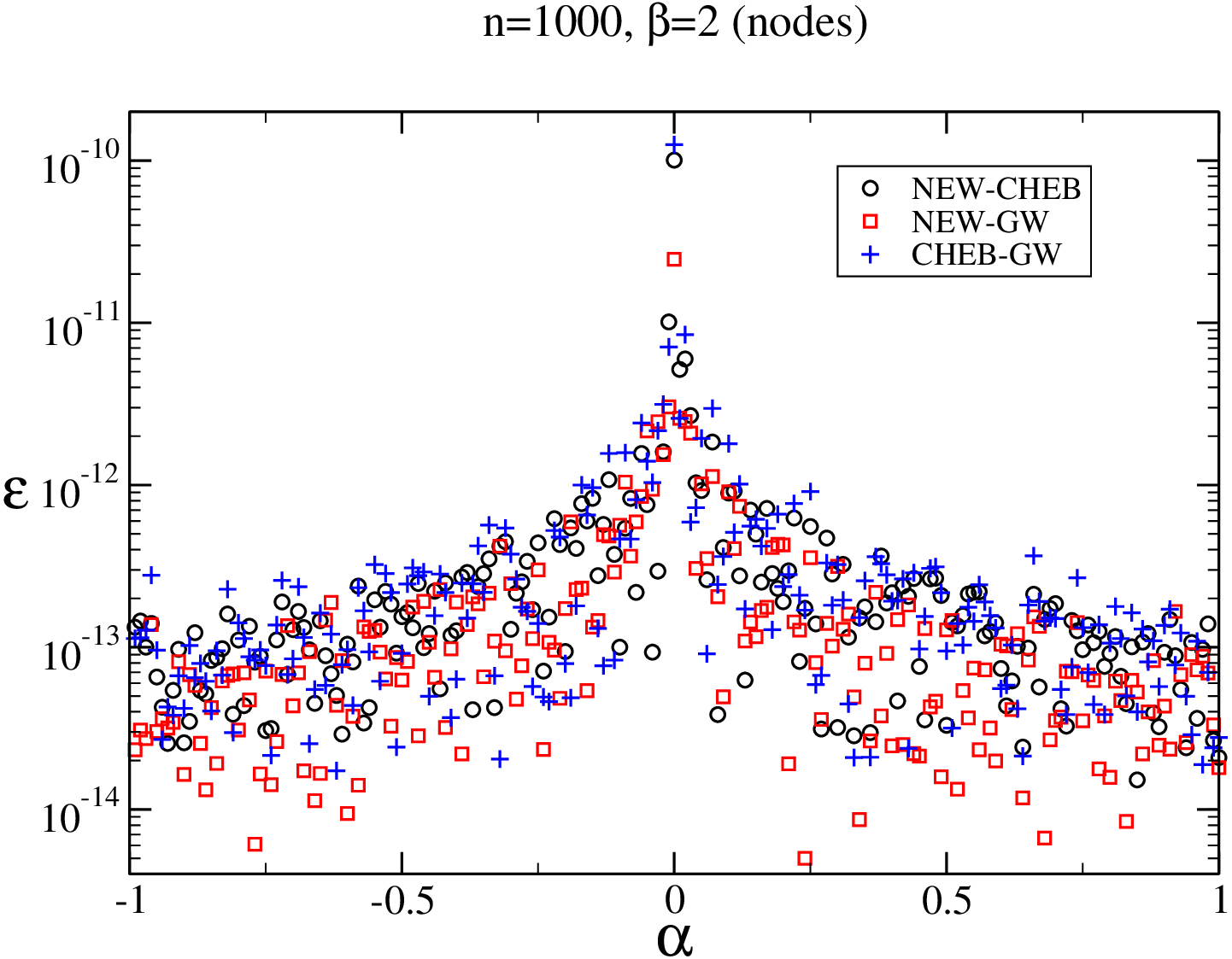}}
\end{minipage}
\end{center}
\end{center}

\caption{Relative maximum accuracies for the weights comparing three different pairs of methods and maximum relative
accuracies for the nodes (bottom, right). In all cases $\beta=2$.}
\label{fig4}
\end{figure}

In Fig. \ref{fig4} we plot the relative maximum accuracy for the weights obtained for Gauss--Jacobi quadrature with parameters
$\alpha\in (-1,1)$ and $\beta=2$ comparing the three different pairs of methods NEW--CHEB, NEW--GW and CHEB--GW, and for three values of $n$:
$n=90$ (when CHEB uses recurrences for computing the polynomials), $n=110$ (when CHEB uses asymptotics) and 
$n=1000$ (again, CHEB uses asymptotics). In addition, we also show the maximum relative accuracy for the nodes in the case $n=1000$ and for
the same values of $\alpha$ and $\beta$. 

For $n=90$, the relative maximum error for the weights is very similar for the three comparisons except that as 
$\alpha\rightarrow -1$ the error for the NEW--GW comparison is smaller; this suggests that 
CHEB loses some accuracy in this limit. The same behavior is observed for the symmetrical cases and we have
checked with our quadruple precision Fortran program that the most accurate method in this limit is NEW.

For $n=110$ CHEB uses asymptotics to compute the polynomials
 and in this case the accuracy of the weights appears to be close to double precision accuracy
 (as checked in the symmetric case by comparing
with our quadruple precision Fortran program). We observe in this case that the error worsens for negative $\alpha$ 
and in particular close to $\alpha=-1/2$; for such values the comparison is favorable for NEW with respect to GW 
and we conclude that 
the most accurate method is CHEB, followed by NEW and the less accurate is GW. We stress again that 
whenever $-1<\alpha,\beta\le 5$ a faster and more accurate method is that of \cite{Gil:2019:NIC}, with close to double precision
accuracy for the nodes and weights. The case $n=1000$ for the weights in Fig. \ref{fig4} shows similar results as for $n=110$.

Finaly, in Fig. \ref{fig4} (bottom, right) 
we show a plot of the maximum relative error in the computation of the nodes when $n=1000$. 
The method NEW is able to compute the nodes with full double precision accuracy and for any $n$; differently, CHEB and GW 
do not compute the nodes with relative accuracy, but with absolute accuracy, which means that there is some relative error 
degradation for the nodes closer to $x=0$ as $n$ increases. For the case shown, this error degradation is more noticeable
for $\lambda$ close to zero because nodes close to $x=0$ occur. We have repeated these tests for other values of $\alpha$ and 
$\beta$ and we conclude that the error degradation for the nodes scales as ${\cal O}(n)$ for CHEB and as ${\cal O}(\sqrt{n})$ 
for GW.

\subsubsection{Comparing against a higher precision algorithm}

In order to confirm the information that we have extracted by intensive comparison tests between three different methods implemented
in double precision (with Matlab), we have also compared the outputs of the NEW and CHEB methods against a high accuracy computation; 
for this purpose, we have implemented the Golub-Welsch algorithm in Maple. These tests are, necessarily, more time 
consuming and less extensive, and the values of $n$ are more limited (for instance the tests for $n=1000$ with Maple are impractical), 
but we use them to illustrate the behaviour of the relative errors for each node and weight. 

\begin{figure}[tb]
\vspace*{0.8cm}
\begin{center}
\begin{minipage}{3cm}
\centerline{\includegraphics[height=4cm,width=6cm]{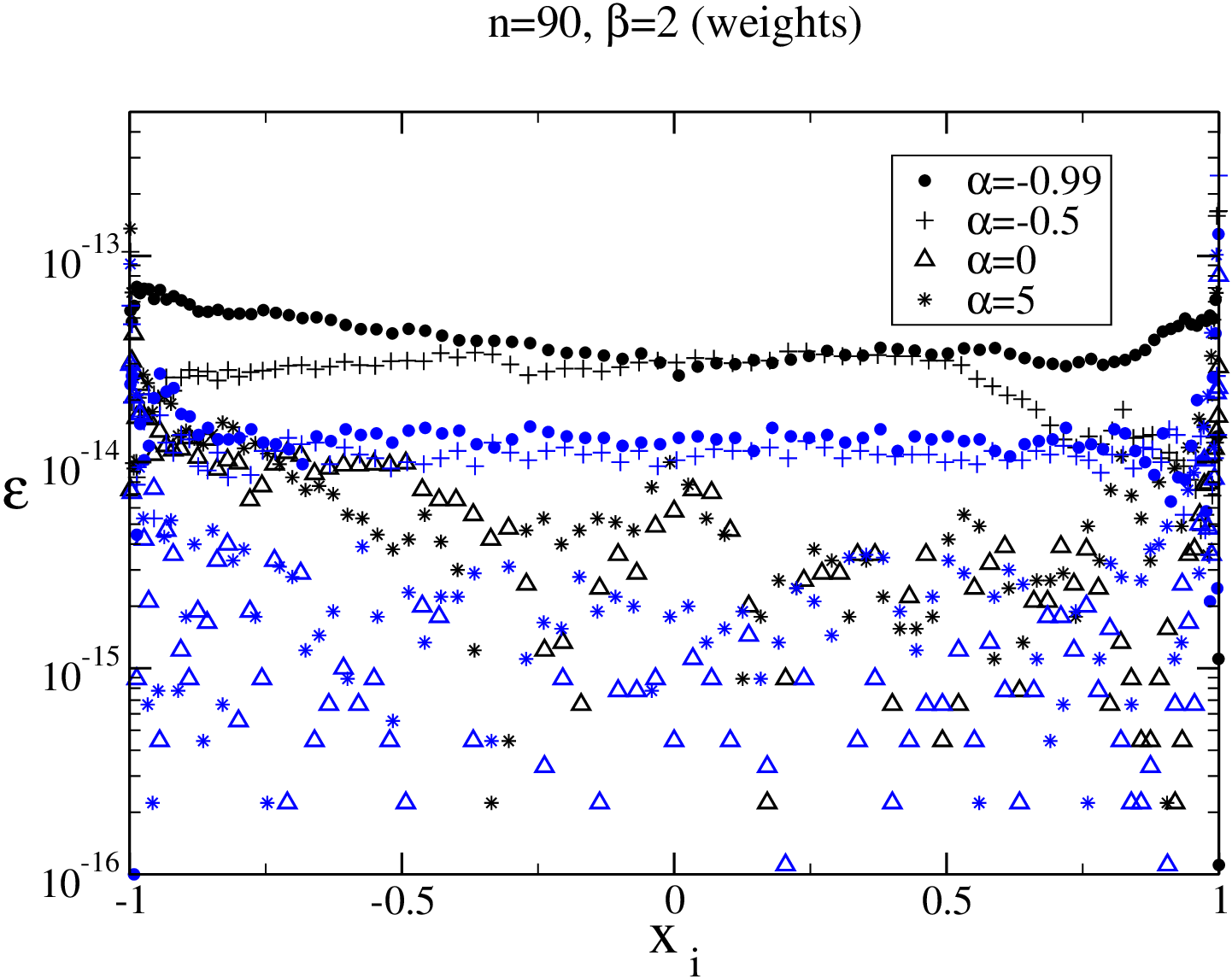}}
\end{minipage}
\hspace*{3.2cm}
\begin{minipage}{3cm}
\centerline{\includegraphics[height=4cm,width=6cm]{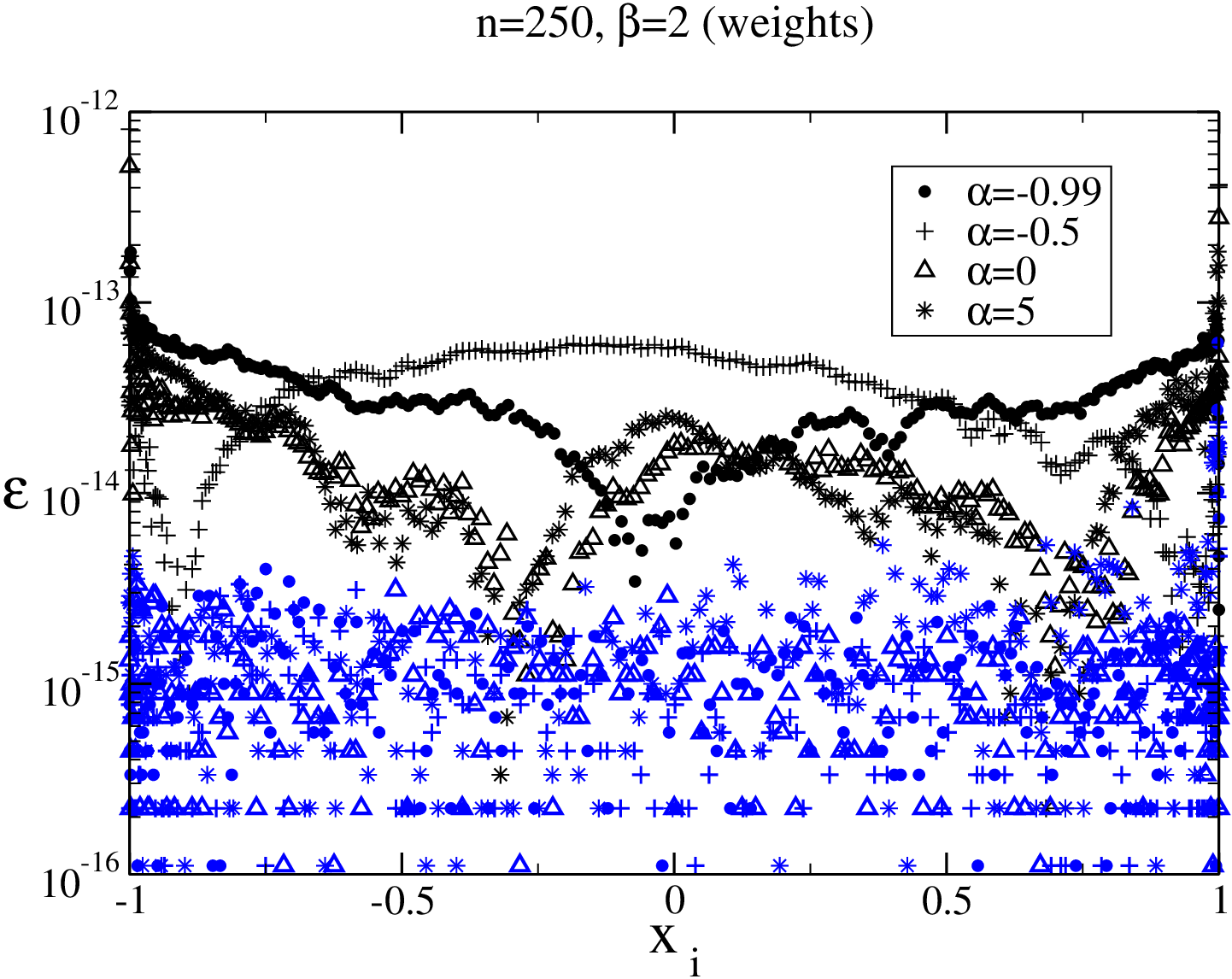}}
\end{minipage}
\end{center}
\vspace*{0.2cm}
\begin{center}
\begin{center}
\begin{minipage}{3cm}
\centerline{\includegraphics[height=4cm,width=6cm]{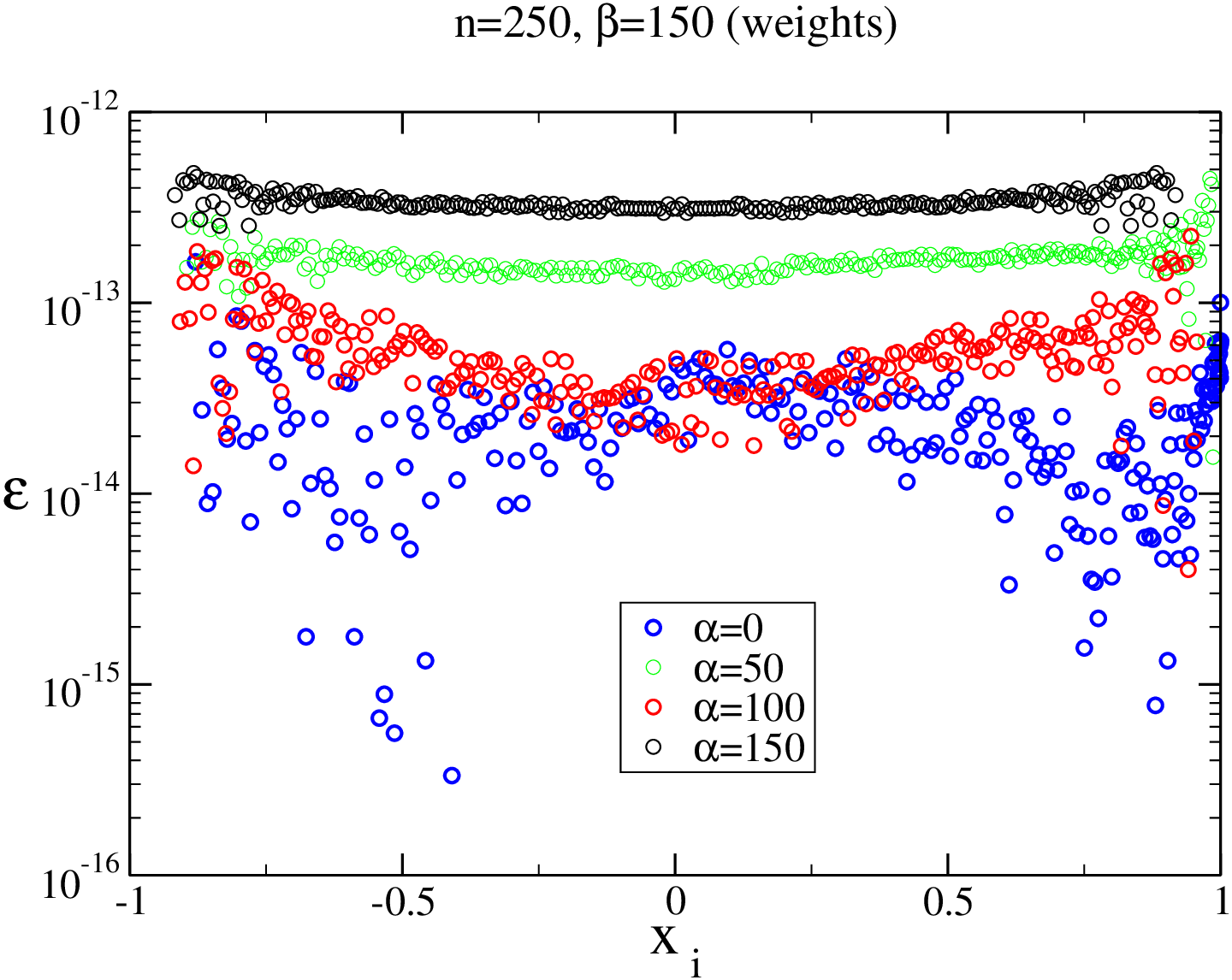}}
\end{minipage}
\hspace*{3.2cm}
\begin{minipage}{3cm}
\centerline{\includegraphics[height=4cm,width=6cm]{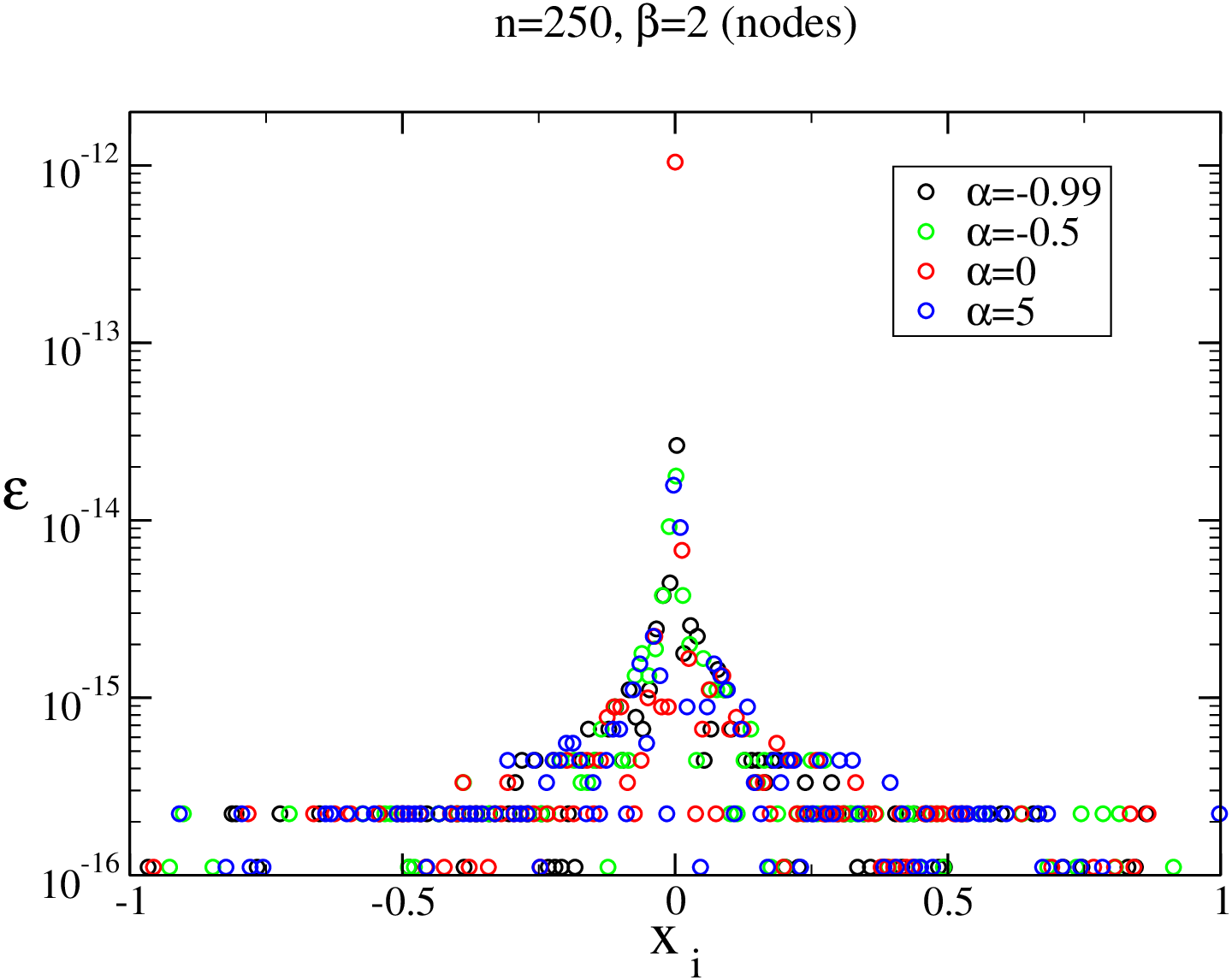}}
\end{minipage}
\end{center}
\end{center}

\caption{Up: relative accuracies for the weights computed with 
the NEW (black) and CHEB (blue) methods for $n=90$ (left) and $n=250$ (right).  Bottom-right: relative accuracy in the computation
 of the nodes with the CHEB algorithm for $n=250$.
 Bottom-left: relative accuracies for the weights computed with the NEW 
method with $n=250$, $\alpha=0, 50, 100, 150$ and $\beta=150$.}
\label{fig5}
\end{figure}

In Fig. \ref{fig5} we show relative errors for the weights. The absolute errors can of course be easily obtained from the relative
errors by multiplying them by the values of the weights (which for large enough $n$ could be estimated with Eq. (\ref{circle})).

The top graphs in Fig. \ref{fig5} display the relative error in the computation of each weight for 
$\beta=2$ and several values of $\alpha$; in the up left figure $n=90$ (when the polynomials are computed by recursion
in CHEB), while $n=250$ in the up right figure (with polynomials computed by asymptotics in CHEB). 
The errors for NEW are ploted in black, while the CHEB errors are plotted in blue.
In the up left figure, we observe 
that CHEB is able to provide typically one additional digit of accuracy with respect to NEW for most the weights; the exception to this is found on the
weights corresponding to the largest nodes (which are the largest weights for $\beta=2$ and $\alpha\le -1/2$, as it is easy to
 check using Eq. (\ref{circle})). 
In particular, we observe a degradation of accuracy as $\alpha\rightarrow -1^{+}$ in CHEB and the relative error for the largest weight
for $n=90$, $\alpha=-0.99$ and $\beta=2$ is $9\times 10^{-11}$ (this point is not shown in the graph in order to improve the visibility 
of the rest of the graph). In the tests corresponding to this up left figure, we used Maple with 40 digits in the computation of
the Goulb-Welsch algorithm.

 In the up right figure, we observe that the relative accuracy for CHEB has improved in relation to NEW with respect to the case
$n=90$,
 and that, again, some accuracy loss happens for the largest weights. For these values of $n>100$, it seems more convenient to use
methods based on asymptotics for computing the weights when they are available as is the case of CHEB (using Newton iterations) 
and also of the methods in \cite{Gil:2019:NIC} (without Newton iterations). 

The bottom-left graph in Fig. \ref{fig5} shows one case for which
asymptotic methods are not available, and for which NEW does produce accurate results. For these large values of the parameters, 
the smallest weights are many orders of magnitude larger that the largest weights, and for this reason  
at least 140 digits of accuracy are needed in Maple in order 
to compute the weights with the Golub-Welsch algorithm, while NEW only requires 15-16 digits.

We summarize in Tables 2-4 the maximum relative and relative maximum errors for the weights for the parameters considering in Fig. 
\ref{fig5}. We observe that when the maximum relative error is equal to the relative maximum errot this means that the relative error
reaches its maximum value for the most significant weight. 

Finally, regarding the errors for the nodes, the absolute errors both fo NEW and CHEB are close to $10^{-16}$ and therefore consistent
with double precision accuracy (the resulting noisy graph is not very interesting and it is not shown). 
The difference between the NEW and CHEB methods is that, as commented before, CHEB computes the nodes
with absolute double precision accuracy but not relative accuracy, which results in some degradation of relative accuracy for the nodes
close to zero. This is shown in Fig. 5 bottom-right. For NEW the relative accuracy is found to be better than $10^{-15}$ and is not 
shown.

As a way of summary, each method has its advantages in terms of accuracy. The main novelty of NEW is that the range of parameters 
available is drastically increased and that, being a method based on finite or convergent processes, it can be extended to 
arbitrary accuracy (as illustrated with the Gauss-Gegenbauer case). Apart from this, as we will see next, the method turns out to be very efficient. 

An optimal algorithm for computing Gauss quadrature rules should combine the
 use of asymptotic methods (when available) with the use of fully convergent methods capable of high accuracy computations, like the one we  have presented. A full and extensive accuracy test of the available methods as a function of the three parameters
 is outside the scope of the present paper.

\begin{table}
$$
\begin{tabular}{rcccccccc}
$n\Rightarrow$  & 90 & 90 & 90 & 90 & 250 & 250 & 250 & 250 \\
$\alpha\Rightarrow$ & -0.99 & -0.5 & 0 & 5 &-0.99  &-0.5 & 0 & 5 \\
\hline
 NEW    & 7.1e-14 & 1.7e-13 & 4.2e-14 & 1.4e-13 &  1.1 e-12 &  8.1e-13 & 5.2e-13 & 1.9e-13\\ 
 CHEB   & 9.6e-11 & 2.5e-13 & 8.1e-14 & 1.0e-13 &  6.1 e-14 &  2.4e-14 & 3.3e-15 & 1.9e-14 
\end{tabular}
$$
\caption{Maximum relative errors $\varepsilon_{mr}$ for the CHEB and NEW methods corresponding to the upper plots of Fig. \ref{fig5}}%
\label{tab2}%
\end{table}

\begin{table}
$$
\begin{tabular}{rcccccccc}
  $n\Rightarrow$  & 90 & 90 & 90 & 90 & 250 & 250 & 250 & 250 \\
  $\alpha\Rightarrow$ & -0.99 & -0.5 & 0 & 5 &-0.99  &-0.5 & 0 & 5 \\
\hline
 NEW    & 3.8e-16 & 1.7e-13 & 3.9e-15 & 5.7e-15 &  2.4e-15 &  4.1e-13 & 1.7e-14 & 1.6e-14\\ 
 CHEB   & 9.6e-11 & 2.5e-13 & 4.4e-15 & 3.1e-15 &  1.4e-15 &  2.4e-14 & 2.1e-15 & 2.7e-15 
\end{tabular}
$$
\caption{Same as table \ref{tab2} but for the relative maximum errors $\varepsilon_{rm}$.}%
\label{tab3}%
\end{table}

\begin{table}
$$
\begin{tabular}{rcccc}
$\alpha\Rightarrow$ & 0 & 50 & 100 & 150\\
\hline
$\varepsilon_{mr}$ & 1.6 e-13 & 4.5 e-13 & 2.2 e-13 & 4.8 e-13 \\
$\varepsilon_{rm}$ & 6.0 e-14 & 1.8 e-13 & 5.0 e-14 & 3.2 e-13 
\end{tabular}
$$
\caption{Relative maximum (absolute) error and maximum relative error for the NEW method with $n=250$, $\beta=150$ and 
$\alpha=0,\,50,\,100,\,150$}
\label{tab3}%
\end{table}

\subsection{CPU times}

\begin{figure}[tb]
\vspace*{0.8cm}
\begin{center}
\begin{minipage}{3cm}
\centerline{\includegraphics[height=4cm,width=6cm]{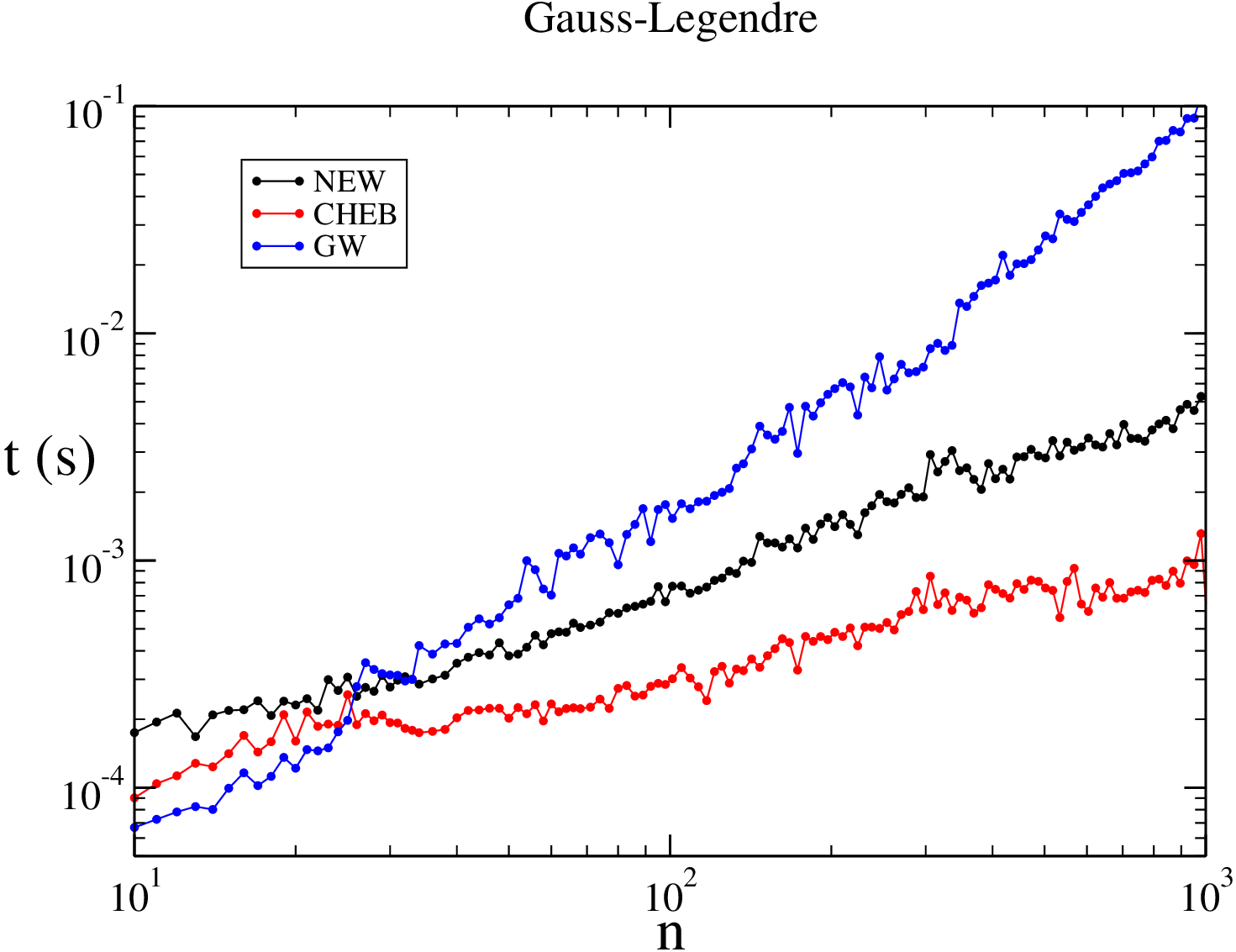}}
\end{minipage}
\hspace*{3.2cm}
\begin{minipage}{3cm}
\centerline{\includegraphics[height=4cm,width=6cm]{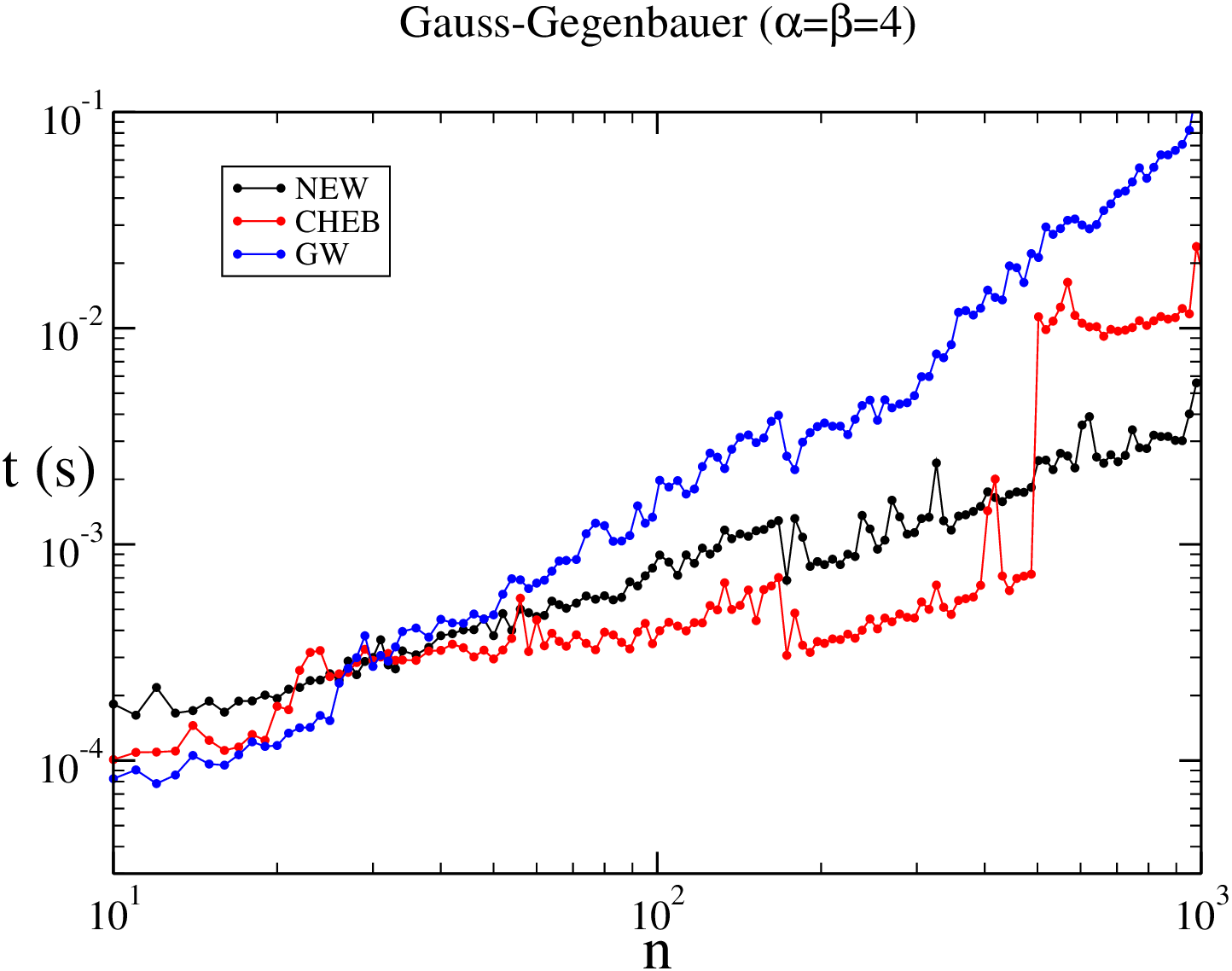}}
\end{minipage}
\end{center}
\vspace*{0.2cm}
\begin{center}
\begin{center}
\begin{minipage}{3cm}
\centerline{\includegraphics[height=4cm,width=6cm]{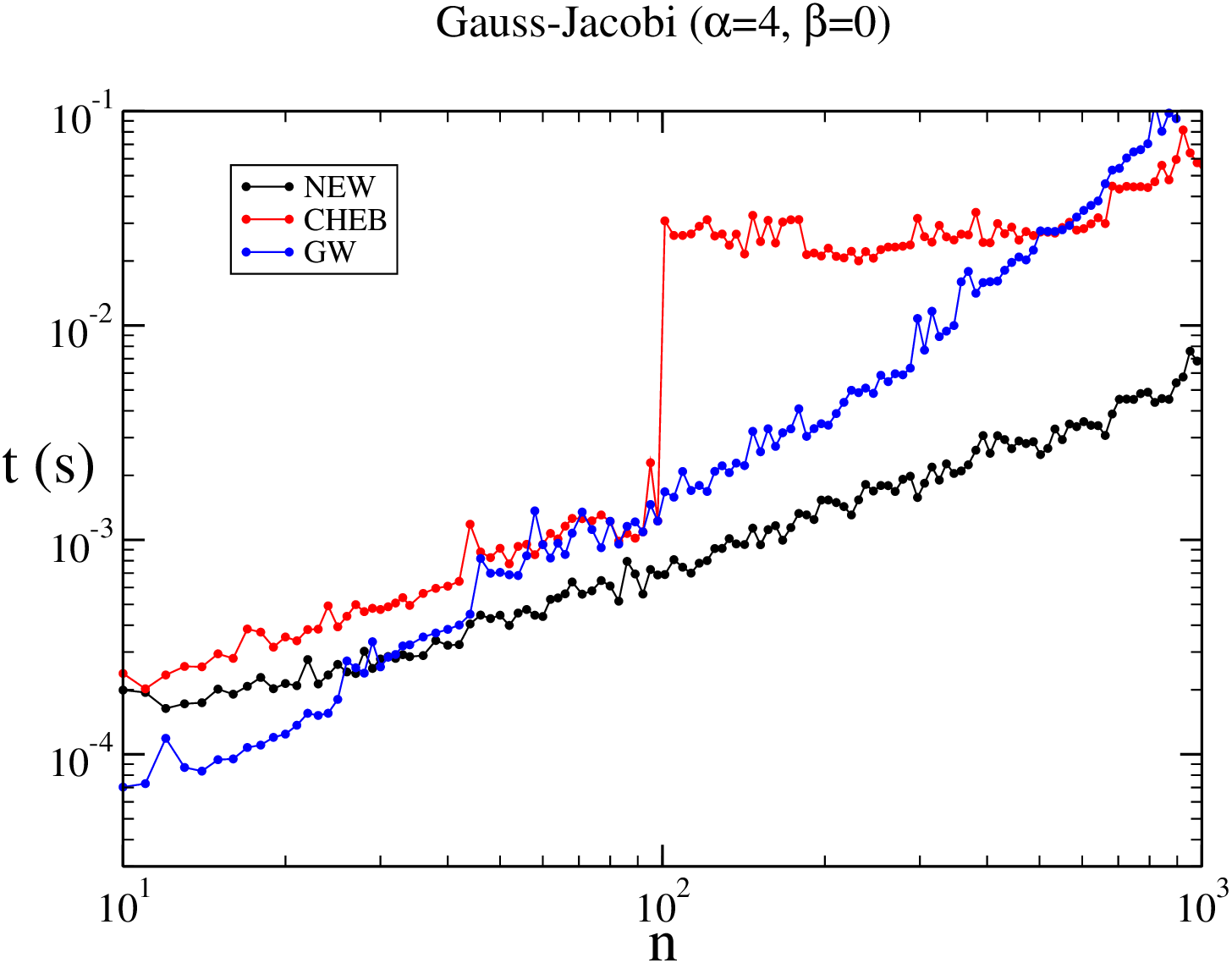}}
\end{minipage}
\hspace*{3.2cm}
\begin{minipage}{3cm}
\centerline{\includegraphics[height=4cm,width=6cm]{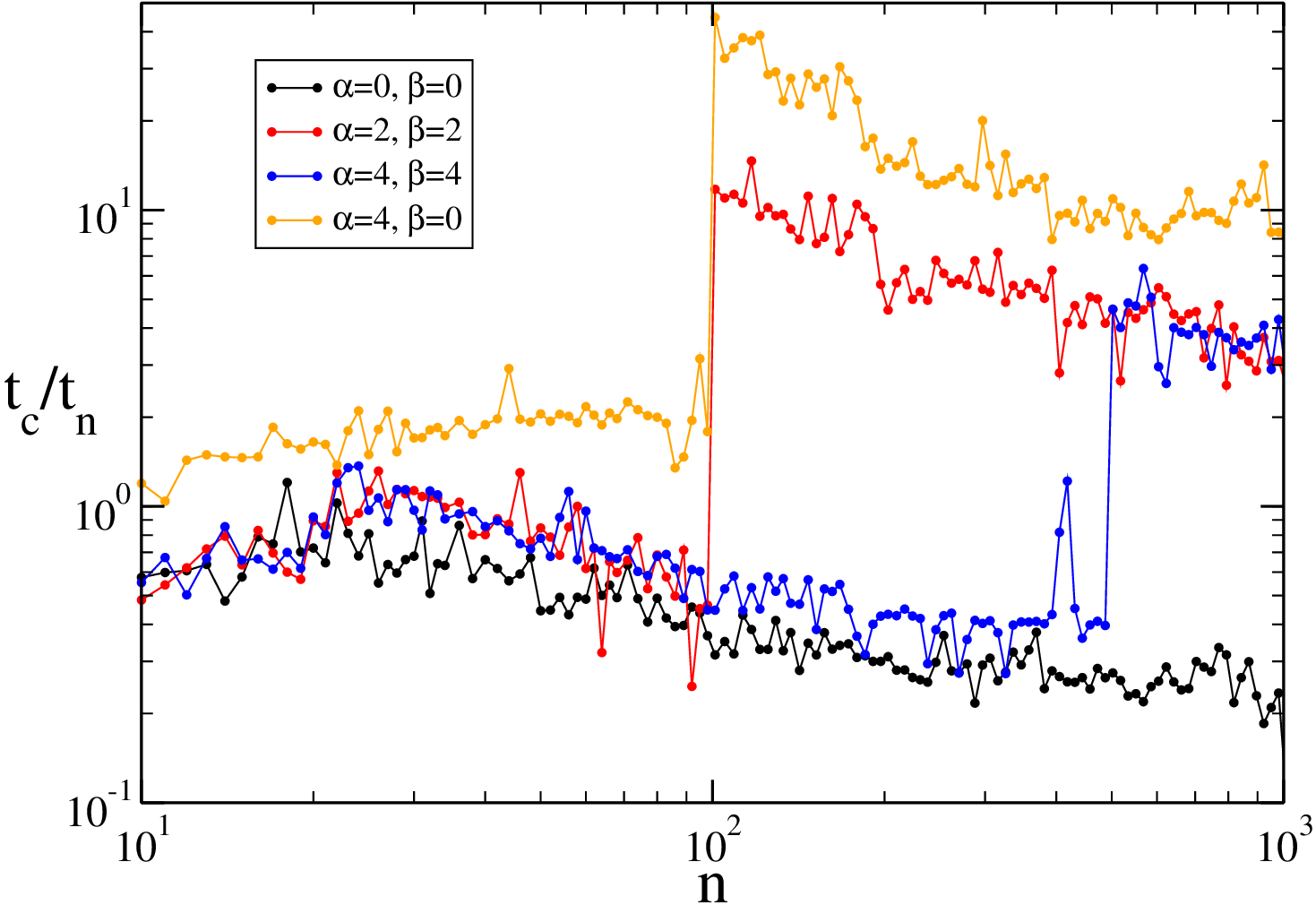}}
\end{minipage}
\end{center}
\end{center}

\caption{CPU time as a function of the degree $n$ for the methods NEW, CHEB and GW and ratio of times between the CHEB and NEW method 
(bottom, right).}
\label{fig6}
\end{figure}

In Fig. \ref{fig6} we compare the CPU time spent by the methods NEW, CHEB and GW as a function of the degree and for several
values of $\alpha$ and $\beta$.

We observe that the behavior of the methods NEW and GW does not change much for the three cases shown (Gauss--Legendre, Gauss--Gegenbauer
 with $\alpha=\beta=4$ and Gauss--Jacobi with $\alpha=4$, $\beta=0$), while for CHEB there are significant differences, because 
CHEB  uses
different methods depending on the value of the parameters and the degree.

For the case of Gauss--Legendre, CHEB uses the asymptotic approximations of \cite{Bogaert:2014:IFC} when $n\ge 100$ and computation
of Legendre polynomials through recurrence for $n<100$; in this case, we observe that CHEB is the fastest method. This could be expected
because for $n>100$ direct asymptotics are used, without iterative methods, and for $n<100$ the simplified expressions for the
Legendre case are also faster to compute. Except for $n<30$, where GW appears to be faster, CHEB appears to be preferable.

For the Gauss--Gegenbauer case shown, CHEB uses the iterative method based on asymptotics for $n>500$ and computation used on recurrences
otherwise; this is observed in the jump in CPU times for $n\ge 500$, and NEW becomes faster in this case. For smaller $n$ the performance
are more or less close to each other. For smaller values of $\lambda$ ($\lambda<3$) the results are quite similar to the next case we
discuss (Gauss--Jacobi), because in this case CHEB uses iteration based on asymptotics for $n\ge 100$ (results not shown).

Finally, for the most general non-symmetric cases, the advantage of NEW in terms of speed becomes clear (Fig. \ref{fig5} bottom, left) 
and only GW is faster for $n<30$. Therefore, except for the symmetric cases, the algorithm NEW is faster. 

Even without making specific algorithms for the symmetric cases in our Matlab implementation, 
it is competitive to CHEB in that cases, and 
clearly faster for non symmetric cases. In order to show this we plot the ratio of the CPU times between the CHEB and NEW methods 
(Fig. \ref{fig5}, bottom, right). We show these ratios, 
as a function of $n$, for four cases: Gauss--Legendre, two Gauss--Gegenbauer cases and 
Gauss--Jacobi example. Our method, except for the Legendre case for large $n$, is competitive when it is not faster.

{\small
\section*{Appendix A: On the conditioning of the recurrence for computing derivatives.}

We briefly discuss the conditioning of the computation of the derivatives with the recurrence relation
(\ref{recun}).
This is a five term recurrence relation, with a space of solutions of dimension four. For studying the
conditioning as $j\rightarrow \infty$ we can divide all terms 
of the recurrence by $j^2$ and then all the coefficients have finite limit as $j\rightarrow +\infty$. We have
$$
\sum_{n=0}^{4} C_n (j) a_{j+2-n}=0,
$$
with
$$
\lim_{j\rightarrow \infty}C_n(j)=c_n=\Frac{Q^{(n)}(x)}{n!}.
$$
Then it is known (see for instance \cite{Elaydi:2005:AIT}, Theorem~8.11) 
that the solutions of the recurrence satisfy
\begin{equation}
\label{lip}
\displaystyle\limsup_{n\rightarrow +\infty}\left(|a_n|\right)^{1/n}=W,
\end{equation}
where $W$ is the modulus of one of the solutions of the characteristic polynomial
$$
\sum_{j=0}^4 c_n \delta^{4-n}=0,
$$
which, upon dividing by $\delta^4$ and denoting $\mu=1/\delta$ we can write
$$
Q(x)+Q'(x)\mu+ Q''(x)\Frac{\mu^2}{2}+Q'''(x)\Frac{\mu^3}{3!}
+Q^{(4)}(x)\Frac{\mu^4}{4!}=0.
$$
And because $Q$ is a polynomial of degree four ($Q(x)=4(1-x^2)^2$) this equation is the same as $Q(x+\mu)=0$,
which, solving for $\mu$ gives two double roots $\mu=x\pm 1$. This means that
the possible values of $W$ in Eq. (\ref{lip}) are $W_1=1/|1-x|$ and $W_2=1/|1+x|$
and there is a subspace of dimension $2$ satisfying (\ref{lip}) with $W=W_1$ and a 
second subspace with $W=W_2$; the first space will be dominant over the second when $W_1>W_2$ and
the opposite when $W_1<W_2$ (the case $W_1=W_2$ is the degenerate case, in which no solution is
exponentially dominant over the rest). 

In our case, we have $a_j=\tilde{Y}^{(j)}(x)$ with $\tilde{Y}(x)$ given by (\ref{Ytil}). We notice that
for $\alpha$ and $\beta$ odd,  $\tilde{Y}(x)$ is a polynomial of degree $N=n+(\alpha+\beta+2)/2=(L+1)/2$, and then
$a_j=0$ if $j>N$. The Taylor series have a finite number of terms in this case and the analysis of stability for $a_j$
as $j\rightarrow \infty$ is not needed. Let us now consider that neither $\alpha$ nor $\beta$ are odd, 
and we leave for later the
case in which only one of the parameters ($\alpha$ or $\beta$) is odd.

When neither $\alpha$ nor $\beta$ are odd, then Taylor series at $x\in (-1,1)$ has an infinite number of terms.
Because $\tilde{Y}(x)$ is a polynomial times $(1-x)^{(\alpha+1)/2}(1+x)^{(\beta +1)/2}$ the radius of convergence
of the series for $\tilde{Y}(x+t)$ centered at $x$,
$$
\tilde{Y}(x+h)=\displaystyle\sum_{j=0}^{\infty} \Frac{\tilde{Y}^{(j)}(x)}{j!} 
h^j=\displaystyle\sum_{j=0}^{\infty} a_j h^j
$$
is $R=\min\{1-x,1+x\}$ (which we could expect because the ODE satisfied by $\tilde{Y}$ has
singularities at $x=\pm 1$) and then
$$
\Frac{1}{R}=\limsup_{n\rightarrow\infty}\sqrt[n]{|a_n|}=\max\{W_1,W_2\}.
$$
Therefore in this case the sequence $\{a_j\}$, $a_j=\tilde{Y}^{(j)}(x)/j!$ is in the dominant subspace
of solutions of the recurrence.

The case when either $\alpha$ or $\beta$ is odd but not both is different. Let us for instance consider that
$\alpha$ is odd, but not $\beta$. In this case, the convergence of the series is limited by the singularity
at $-1$ (but not at $+1$); the radius of convergence in this case is therefore $R=1+x$ and then
$$
\limsup_{n\rightarrow\infty}\sqrt[n]{|a_n|}=\Frac{1}{|1+x|}.
$$ 
Therefore, $\{a_j\}$ is in the dominant subspace only if $x\in (-1,0)$. In this case for large enough $j$
the forward computation of $a_j$ would be unstable for positive $x$. 
However, even for this case we have not observed
inaccuracies in the computation of the series. For a given accuracy claim, the number of terms in the
series needed are not high enough to produce stability issues.  
}

\begin{acknowledgements}
NMT thanks CWI  for scientific support. 
\end{acknowledgements}

\bibliographystyle{spmpsci}
\bibliography{gaussr}

\end{document}